\documentclass[11pt, a4paper]{article}
\usepackage{latexsym,amssymb}
\usepackage{amsmath}
\usepackage[mathscr]{eucal}
\usepackage{color, xcolor}
\usepackage[abbrev]{amsrefs}

\newtheorem{Theorem}{\bf Theorem}[section]
\newtheorem{Lemma}{\bf Lemma}[section]
\newtheorem{Proposition}{\bf Proposition}[section]
\newtheorem{Corollary}{\bf Corollary}[section]
\newtheorem{Remark}{\bf Remark}[section]
\newtheorem{Example}{\bf Example}[section]
\newtheorem{Definition}{\bf Definition}[section]

\newenvironment{theorem}{\begin{Theorem}$\!\!\!$}{\end{Theorem}}
\newenvironment{lemma}{\begin{Lemma}$\!\!\!$}{\end{Lemma}}
\newenvironment{proposition}{\begin{Proposition}$\!\!\!$}{\end{Proposition}}

\newenvironment{remark}{\begin{Remark}$\!\!\!$}{\end{Remark}}

\newenvironment{definition}{\begin{Definition}$\!\!\!$}{\end{Definition}}

\numberwithin{equation}{section}

\numberwithin{equation}{section}

\newcommand{\dee}{{\rm{d}}}

\newcommand{\R}{\mathbb{R}}
\newcommand{\N}{\mathbb{N}}
\newcommand{\E}{\varepsilon}
\newcommand*{\RL}{\partial_{\mathnormal t}^{\alpha}}

\newcommand*{\dt}{\dfrac{\dee}{\dee{\mathnormal t}}}

\DeclareMathOperator*{\esssup}{ess\sup}
\DeclareMathOperator*{\essinf}{ess\inf}

\begin{document}
\title{Existence and uniqueness of $L^1$-solutions to\\ time-fractional nonlinear diffusion equations}
\author{Mikiya Kametaka and Tatsuki Kawakami}
\date{}
\maketitle
\begin{abstract}
We establish the global existence and uniqueness of $L^1$-solutions to
the Cauchy problem for time-fractional porous medium type nonlinear diffusion equations.
Furthermore, we give the mass conservation law for $L^1$-solutions
to time-fractional fast diffusion equations,
and prove that the finite-time extinction does not occur for any nonnegative $L^1$-solutions.
\end{abstract}
\vspace{10pt}
\noindent
Keywords: 
time-fractional diffusion, porous medium equation, fast diffusion equation, \\
$L^1$-solution, $L^1$-contraction principle, mass conservation law, nonextinction in finite time

\vspace{40pt}
\noindent 
Addresses:

\noindent 
M. K.: Graduate School of Applied Mathematics and Informatics Course,\\ 
Faculty of Advanced Science and Technology, Ryukoku University,\\
1-5 Yokotani, Seta Oe-cho, Otsu, Shiga 520-2194, Japan\\
\noindent 
E-mail: {\tt Y25D001@mail.ryukoku.ac.jp}\\

\smallskip
\noindent 
{T. K.}: Applied Mathematics and Informatics Course,\\ 
Faculty of Advanced Science and Technology, Ryukoku University,\\
1-5 Yokotani, Seta Oe-cho, Otsu, Shiga 520-2194, Japan\\
\noindent 
E-mail: {\tt kawakami@math.ryukoku.ac.jp}\\

\newpage

\section{Introduction}
Consider the Cauchy problem for time-fractional porous medium type nonlinear diffusion equations
\begin{equation}
\label{eq:P}\tag{P}
	\begin{cases}
		\partial_t\left(k*[u-u_0]\right) = \Delta\Phi(u), & x \in \R^N,\ \ t>0, 
		\vspace{5pt}\\
		u(x,0) = u_0(x), & x \in \R^N,
	\end{cases}
\end{equation}
where $N\ge1$, $\partial_t := \partial/\partial t$, 
and $u_0$ is an integrable function in $\R^N$.
Here $k*v$ denotes the convolution on the positive halfline with respect to the time variable,
that is,
$$
	(k*v)(t):=\int_0^tk(t-\tau)v(\tau)\,\dee \tau\qquad\mbox{for $k,v\in L^1_{\rm loc}([0,\infty))$},\qquad t\ge0.
$$
Throughout this paper, we assume that
\begin{enumerate}
\item[({$K$})]
	$k \in L^1_{\rm loc}([0,\infty))$ is a nonnegative nonincreasing function,
	and there exists a nonnegative nonincreasing function $\ell \in L^1_{\rm loc}([0,\infty))$ such that
	$(k*\ell)(t)\equiv1$ for $t\ge0$.
\end{enumerate}
Hence $k$ is a completely positive kernel (see, e.g., \cite{C}*{Theorem~2.2}).
A typical example of $k$ satisfying the assumption ($K$) is so-called the Riemann-Liouville kernel, that is,
$$
	k(t) = g_{\alpha}(t):=\frac{t^{-\alpha}}{\Gamma(1-\alpha)},\quad 0<\alpha<1.
$$
In this case, the nonlocal time-derivative $\RL u:=\partial_t(g_\alpha*u)$ 
is so-called the $\alpha$-th order Riemann-Liouville derivative of $u$,
and $\RL (u-u_0)$ is so-called the $\alpha$-th order Caputo derivative of $u$ if it is sufficiently smooth.
Furthermore, for the function $\Phi$, we assume that
\begin{enumerate}
\item[({$L$})]
	$\Phi \in C(\R) \cap C^1(\R\setminus\{0\})$ is symmetric with respect to the origin, 
	namely, $\Phi(r)=-\Phi(-r)$ for $r>0$, and strictly increasing on $\R$ with $\Phi(0)=0$.
\end{enumerate}
A typical example of $\Phi$ satisfying the assumption ($L$) is given by
$$
	\Phi(r) = |r|^{m-1}r,\quad m>0.
$$
In this case, the equation \eqref{eq:P} is so-called the fast diffusion equation if $0<m<1$, 
the heat equation if $m=1$, and the porous medium equation if $m>1$.
In this paper, under conditions $(K)$ and $(L)$,
we study the $L^1$-theory,
in particular, the existence and uniqueness of $L^1$-solutions to problem \eqref{eq:P}, which is defined as follows.
\begin{definition}
\label{Definition:1.1}
	Let $u_0\in L^1(\R^N)$.
	We say that $u$ is an $L^1$-solution to problem \eqref{eq:P} if
	\begin{itemize}
		\item[\rm (i)] 
		$u \in L^\infty(0,\infty;L^1(\R^N))$,
		\item[\rm (ii)] 
		$\Phi(u) \in L_{\rm loc}^1(0,\infty; L_{\rm loc}^1(\R^N))$ with 
		$\ell*\|\Phi(u)\|_{L_{\rm loc}^1(\R^N)} \in L^\infty(0,\infty)$,
		\item[\rm (iii)] 
		$u$ satisfies the identity
		\begin{equation}
		\label{eq:1.1}
			\int_{\R^N}\left(u(x,t)-u_0(x)\right)\phi(x)\,\dee x
			=\left(\ell*\left[\int_{\R^N}\Phi(u(x,\cdot))\Delta\phi(x)\,\dee x\right]\right)(t),\,\,{\rm a.e.}\,\,t>0,
		\end{equation}
		for any test function $\phi\in C_0^\infty(\R^N)$.
	\end{itemize}
\end{definition}
\medskip

For classical porous medium type nonlinear diffusion equations
\begin{equation}
\label{eq:PME}\tag{PME}
	\partial_tu = \Delta\Phi(u),
\end{equation}
it has been given extensive studies such as solvability, regularity, and so on.
In particular, the $L^1$-theory for \eqref{eq:PME} is one of the most fundamental and important issues, 
and it has been well established in many papers (e.g., \cites{BC, CP, HP} and the references therein) 
and in survey books \cites{V1,V2}.
It is well known that, for any $u_0 \in L^1(\R^N)$,
there exists a global-in-time solution $u$ of \eqref{eq:PME}, which belongs to $C([0,\infty);L^1(\R^N))$.
Furthermore, for two solutions $u,v$ of \eqref{eq:PME},  it holds that
\begin{equation}
\label{eq:1.2}
	\int_{\R^N}[u(x,t)-v(x,t)]_+\,\dee x \le \int_{\R^N}[u(x,s)-v(x,s)]_+\,\dee x,\ \ t \ge s \ge 0,
\end{equation}
where
$$
	[r]_+=
	\begin{cases}
	r & {\rm if}\,\,r > 0,\vspace{5pt}\\
	0 & {\rm if}\,\,r \le 0.
	\end{cases}
$$
This estimate is so-called the {\it $L^1$-contraction principle},
and applying \eqref{eq:1.2}, we can obtain uniqueness of solutions and the comparison principle.
Moreover, for the case $\Phi(r)=|r|^{m-1}r$ with $m>0$, the mass conservation law, namely,
$$
	\int_{\R^N}u(x,t)\,\dee x = \int_{\R^N}u_0(x)\,\dee x,\qquad t>0,
$$
holds for any $u_0 \in L^1(\R^N)$ if the exponent $m$ satisfies that
$$
	m > 0 \quad{\rm if}\quad N=1,2,\qquad
	m \ge \frac{N-2}{N}\quad {\rm if}\quad N\ge3.
$$
See, e.g., \cite{V}. 

On the other hand, time-fractional diffusion equations have been recently studied by many authors 
(e.g., \cites{A,AN,AAC,AB,BGI,CQW,EK,KSZ,KSVZ,KY,L,LY,PL,SW,VWZ,VZ1,VZ2,VZ3,WWZ,Z1,Z2,Z3,Z4}).
Among others, for the case where $\Omega$ is bounded, 
Wittbold, Wolejko, and Zacher \cite{WWZ} proved 
the existence of bounded weak solutions to time-fractional degenerate quasilinear diffusion equations
with measurable coefficients
$$
	\partial_t(k*[u-u_0])=\mbox{div}(A(x,t)\nabla\varphi(u))
$$
under the homogeneous Dirichlet boundary condition (see Definition~\ref{Definition:3.1} and Proposition~\ref{Proposition:3.1}).
They also showed that, under an additional regularity assumption on solutions, namely
$$
	k*[u-u_0] \in W^{1,1}(0,T; L^1(\Omega))\,\,\,{\rm with}\,\,\,(k*[u-u_0])(0)=0,
$$
the local $L^1$-contraction type principle
$$
	\|u-v\|_{L^1(0,T;L^1(\Omega))} \le T\|u_0-v_0\|_{L^1(\Omega)}
$$
holds for each $T>0$.
However, the above space does not coincide with the one of the weak solutions they constructed,
and the existence of solutions belonging to such a space has not been established.  
Furthermore, Akagi \cite{A} established the abstract theory for time-fractional gradient flow 
of subdifferential type in a real Hilbert space $H$
$$
	\dt\left(k*[u-u_0]\right)(t) + \partial\varphi(u(t)) \ni 0\quad{\rm in}\quad H,\quad t>0,
$$
and proved the existence and uniqueness of strong solutions to the above equation (see also \cites{AN,LL}).
Here $\partial\varphi$ denotes the subdifferential operator of $\varphi$, 
which is applicable to porous medium and fast diffusion equations,
the $p$-Laplace parabolic equation, the Allen--Cahn equation, and so on.
Recently, Bonforte, Gualdani, and Ibarrondo \cite{BGI} extended this result to 
nonlocal porous medium type equation
$$
	\RL(u-u_0) = -\mathcal{L}(|u|^{m-1}u),\quad m > 0.
$$
A typical example of nonlocal operators $\mathcal{L}$ is a fractional Laplacian $(-\Delta)^s$ with $s\in(0,1)$.
They prove the existence and uniqueness of $H^*$-solutions (see \cite{BGI}*{Definition~1.1}, 
and this notion coincides with weak solutions defined by Definition~\ref{Definition:3.1} if $\mathcal{L}=-\Delta$), 
the comparison principle and $L^p-L^\infty$ smoothing effects with $p>1$.
Furthermore, they also prove that the finite-time extinction does not occur for any nonnegative $H^*$-solutions
if $0<m<1$.
\medskip

To the best of our knowledge,
there are almost no results on the $L^1$-theory for time-fractional nonlinear diffusion equations,
in particular for the case of the whole space $\R^N$.
In this paper, we establish the existence and uniqueness of $L^1$-solutions to problem~\eqref{eq:P},
and we also prove the comparison principle.
Furthermore, we also show that,
for the case $\Phi(r)=|r|^{m-1}r$ with $0<m<1$,
the mass conservation law holds for any initial data $u_0 \in L^1(\R^N)$.
Moreover, for the time-fractional fast diffusion equation
\begin{equation}
\label{eq:TFDE}
\tag{TFDE}
	\RL(u-u_0)=\Delta u^m,\quad 0<m<1,
\end{equation}
we prove that the finite-time extinction does not occur for any nonnegative $L^1$-solutions.
\medskip

We first give the existence and uniqueness of $L^1$-solutions to problem~\eqref{eq:P},
which satisfy the $L^1$-contraction principle.
\begin{theorem}
\label{Theorem:1.1}
	Assume $(K)$ and $(L)$.
	Let $u_0 \in L^1(\R^N)$.
	\begin{itemize}
		\item[\rm (a)]
		There exists a unique $L^1$-solution $u\in C([0,\infty);L^1(\R^N))$ to problem \eqref{eq:P}
		satisfying
		$$
			\|u\|_{L^\infty(0,\infty; L^1(\R^N))} \leq \|u_0\|_{L^1(\R^N)}.
		$$
		Furthermore, it holds that
		$$
			\ell*\|\Phi(u)\|_{L_{\rm loc}^1(\R^N)} \in C([0,\infty)) \cap L^\infty(0,\infty)\,\,\,{\rm with}\,\,\,
			\left(\ell*\|\Phi(u)\|_{L_{\rm loc}^1(\R^N)}\right)(0)=0.
		$$
		\item[\rm (b)]
		Let $u,v \in C([0,\infty);L^1(\R^N))$ be $L^1$-solutions to problem \eqref{eq:P} 
		with initial data $u_0,v_0 \in L^1(\R^N)$, respectively.
		Then it holds that
		\begin{equation}
		\label{eq:1.3}
			\int_{\R^N}[u(x,t)-v(x,t)]_+\,\dee x \le \int_{\R^N}[u_0(x)-v_0(x)]_+\,\dee x,\quad t\ge0.
		\end{equation}
		In particular,
		\begin{equation}
		\label{eq:1.4}
			\|u(t)-v(t)\|_{L^1(\R^N)} \le \|u_0-v_0\|_{L^1(\R^N)},\quad t \ge 0.
		\end{equation}
	\end{itemize}
\end{theorem}
\medskip
\begin{remark}
\label{Remark:1.1}
	{\rm (i)}
	For the case where $\Omega$ is bounded,
	if $u_0\in L^\infty$, then, applying improved energy estimates 
	given in Theorem~$\ref{Theorem:3.1}$ with Poincar\'e's inequality, 
	we see that weak solutions provided in \cites{A, BGI, WWZ} 
	$($see also Definition~$\ref{Definition:3.1}$$)$,
	become $L^1$-solutions as defined in Definition~$\ref{Definition:3.2}$.
	\newline
	{\rm (ii)}
	By assertion~$($b$)$ of Theorem~$\ref{Theorem:1.1}$,
	we see that
	the comparison principle for $L^1$-solutions to \eqref{eq:P} holds,
	that is, if $u_0 \le v_0$ a.e. $x\in\R^N$, then $u \le v$ a.e. $x\in\R^N$ and $t>0$.
	In particular, if $u_0 \ge 0$, then $u \ge 0$ a.e. $x\in\R^N$ and $t>0$.
\end{remark}

The proof is based on the idea of \cite{V2}*{Theorem~9.3}.
More precisely, we begin by constructing $L^1$-solutions in a ball $B_R:=\{x\in\R^N:|x|<R\}$,
and then, by taking the limit as $R\to\infty$, we obtain desired solutions.
To this end, we first improve energy estimates for weak solutions given in \cite{WWZ}.
Next we establish the $L^1$-theory in a bounded domain by using above solutions as approximate ones.
In particular, we show the existence and uniqueness of $L^1$-solutions satisfying a weighted $L^1$-contraction principle.
Finally, taking suitable weight functions and the limit as $R\to\infty$, we prove Theorem~\ref{Theorem:1.1}.
\medskip

From now we focus on the fast diffusion case, that is, $\Phi(r)=|r|^{m-1}r$ with $m \in (0,1)$,
and show several properties of $L^1$-solutions.
We next show that the mass conservation law holds under the same condition on $m$
as in \eqref{eq:PME}, which involves {\it local} time derivative.
\begin{theorem}
\label{Theorem:1.2}
	Let $\Phi(r)=|r|^{m-1}r$ with
	$$
		0<m<1 \quad {\rm if}\quad N=1,2,\qquad
		\frac{N-2}{N} \le m<1\quad{\rm if}\quad N\ge3.
	$$
	Assume $(K)$. 
	Let $u \in C([0,\infty);L^1(\R^N))$ be an $L^1$-solution to problem \eqref{eq:P} 
	with initial data $u_0\in L^1(\R^N)$.
	Then it holds that
	\begin{equation}
	\label{eq:1.5}
		\int_{\R^N}u(x,t)\,\dee x=\int_{\R^N}u_0(x)\,\dee x,\quad t>0.
	\end{equation}
\end{theorem}
\begin{remark}
\label{Remark:1.2}
	$\rm(i)$
	For $m>1$, namely the porous medium case,
	 if initial data $u_0$ also belongs to $L^p(\R^N)$ with $p \ge m$, 
	 then the mass conservation law \eqref{eq:1.5} holds 
	 $($see below Remark~$\ref{Remark:4.1}$ for details$)$.
	 \newline
	$\rm(ii)$
	For the case of local diffusion with the Caputo derivative, that is, the time-fractional heat equation
	$$
	\RL(u-u_0)=\Delta u,
	$$
	the fundamental solution can be computed explicitly.
	Therefore, applying the explicit representation formula of solutions using this fundamental solution,
	we can easily show that the mass conservation law \eqref{eq:1.5} holds for $u_0\in L^1(\R^N)$
	$($see \cites{CQW, KSZ}$)$.
\end{remark}

Finally, for the case of the Caputo derivative, namely \eqref{eq:TFDE},
we show that nonnegative $L^1$-solutions do not extinguish in finite time. 
\begin{theorem}
\label{Theorem:1.3}
	Let $p\in[1,\infty]$, and let
	$u \in C([0,\infty); L^1(\R^N)) \cap L^\infty(0,\infty;L^p(\R^N))$ be a nonnegative $L^1$-solution 
	of \eqref{eq:TFDE} with nonnegative initial data $u_0 \in L^1(\R^N) \cap L^p(\R^N)$.
	Then it holds that
	\begin{equation}
	\label{eq:1.6}
		\|u(t)\|_{L^p(\R^N)}\ge C\left(1+t^{\frac{\alpha}{m}}\right)^{-1},\quad t>0,
	\end{equation}
	where $C$ depends only on $m$, $N$, $\alpha$, $p$, and $\|u_0\|_{L^p(\R^N)}$.
\end{theorem}
\begin{remark}
\label{Remark:1.3}
	$\rm(i)$
	It is well known that solutions of \eqref{eq:PME}, namely $\alpha=1$, extinguish in finite time, 
	that is, there exists a $T_0>0$ such that
	$$
		\|u(t)\|_{L^\infty(\R^N)} \equiv 0,\quad t\ge T_0
	$$
	if initial data $u_0$ belongs to $L^{p_*}(\R^N)$ with $p_*=N(1-m)/2>1$ and $N\ge3$, that is,
	$0<m<(N-2)/N$ $($see \cite{V1}$)$.
	\vspace{3pt}
	\newline
	$\rm(ii)$
	For the case where $\Omega$ is bounded, in \cite{BGI}, 
	Bonforte {\it et al.} proved that a similar estimate to \eqref{eq:1.6} holds 
	for nonnegative weak solutions of \eqref{eq:TFDE}, that is,
	$$
		\|u(t)\omega_1\|_{L^1(\Omega)}\ge C\left(1+t^{\frac{\alpha}{m}}\right)^{-1},\quad t>0,
	$$
	where $\omega_1$ is the first eigenfunction of $(-\Delta)_\Omega$ $($see \cite{BGI}*{Proposition 5.2}$)$.
\end{remark}

By Theorem~\ref{Theorem:1.3} and Remark~\ref{Remark:1.3},
we see that the extinction of solutions is not continuous at $\alpha=1$,
and this is one of the typical effects of the time nonlocality, often referred to as the ``memory effect".
In fact, this discontinuity is caused by the behavior of solutions to the fractional ODE
$$
	\RL(v-v_0)=-\lambda v^m,\quad v(0)=v_0>0,
$$
where $\lambda>0$.
If $\alpha=1$ and $m\in(0,1)$, then we can easily see that
$$
	v(t)=\biggl[v_0^{1-m}-\lambda(1-m)t\biggr]_+^{\frac{1}{1-m}},\quad t>0,
$$
and this implies the finite time extinction.
On the other hand,
if $0<\alpha<1$ and $m>0$, then there exists a constant $C>0$ such that
$$
	C^{-1}\left(1+t^{\frac{\alpha}{m}}\right)^{-1} \le v(t) 
	\le C\left(1+t^{\frac{\alpha}{m}}\right)^{-1},\quad t>0
$$
(see \cite{VZ2}*{Theorem~7.1}).
This lower estimate plays a crucial role to prove Theorem~\ref{Theorem:1.3}.
\medskip

The rest of this paper is organized as follows.
In Section \ref{section:2}, we recall basic facts on nonlocal time-derivative operators.
In Section \ref{section:3}, 
we consider the initial-boundary value problem in a bounded domain with Dirichlet zero condition.
We first establish the existence and uniqueness of weak solutions 
for degenerate and singular problems.
Moreover, we develop the $L^1$-theory for bounded domains.
In particular, we construct $L^1$-solutions and provide the weighted $L^1$-contraction principle.
Section \ref{section:4} is devoted to proofs of main results: 
we first establish the existence and uniqueness of $L^1$-solutions to problem \eqref{eq:P}, 
and prove the mass conservation law for $\Phi(r)=|r|^{m-1}r$, $0<m<1$.
We finally show the nonextinction property of nonnegative $L^1$-solutions 
to the time-fractional fast diffusion equation \eqref{eq:TFDE}.

	\section{Preliminaries}\label{section:2}
In this section, we recall basic facts on nonlocal time-derivative operators.

We first recall some properties with respect to the Yosida approximation.
Let $T>0$, $1\le p<\infty$, and let $X$ be a real Banach space.
We define the operator $B$ by
$$
	Bu = \dt\left(k*u\right),\quad 
	D(B)=\left\{u \in L^p(0,T; X); k*u \in W^{1,p}(0,T; X)\,\,\,\mbox{with}\,\,\,(k*u)(0)=0\right\}.
$$
It is well known that $B$ is $m$-accretive operator in $L^p(0,T; X)$ (see \cites{B,CP,GLS,G}).
Let $n \in \N$ and $B_n$ denotes the Yosida approximation of $B$, which is defined by $B_n = nB(n+B)^{-1}$.
Then, for any $u \in D(B)$, it holds that
\begin{equation}
\label{eq:2.1}
	B_nu \to Bu\quad\mbox{strongly in $L^p(0,T; X)$}\quad\mbox{as}\quad n \to \infty.
\end{equation}
Moreover, $B_n$ admits the representation
$$
	B_nu=\dt\left(k_n*u\right),\quad u \in L^p(0,T; X),
$$
where $k_{n} = ns_{n}$, and $s_{n}$ is the unique solution to the Cauchy problem
$$
		\dt\left(k*[s_n-s_n(0)]\right) + ns_n=0,\quad t > 0,\qquad
		s_n(0) = 1
$$
(see \cite{VZ1}).
This problem can also be written as the Volterra integral equation
\begin{equation}
\label{eq:2.2}
	s_{n}(t) + n(s_{n}*\ell)(t) =1,\quad t>0.
\end{equation}
By the standard theorem to the Volterra integral equation (see, e.g., \cite{P}*{Proposition~4.5}),
we see that $s_n \in W_{\rm loc}^{1,1}([0,\infty))$ is a nonnegative nonincreasing function.
Let $h_n \in L_{\rm loc}^1([0,\infty))$ be the nonnegative resolvent kernel associated with $n\ell$, that is,
\begin{equation}
\label{eq:2.3}
	h_{n}(t) + n(h_{n}*\ell)(t) = n\ell(t),\quad t\in(0,T)
\end{equation}
(see, e.g., \cite{GLS}).
Then it follows from $k*\ell=1$ that
$$
	(k*h_n)(t) + n(k*h_n*\ell)(t) = n,\quad t\in(0,T),
$$
which together with \eqref{eq:2.2} implies that
\begin{equation}
\label{eq:2.4}
	k_n = ns_n = k*h_n.
\end{equation}
Note that, for any function $f \in L^p(0,T; X)$, it holds that
\begin{equation}
\label{eq:2.5}
	h_n*f \to f\quad\mbox{strongly in $L^p(0,T; X)$}\quad\mbox{as}\quad n \to \infty.
\end{equation}
In fact, set $v:=\ell*f$. 
Then, by \eqref{eq:2.1}, we see that $v \in D(B)$ and
$$
	\dt(k_n*v)=B_nv 
	\xrightarrow[n\to\infty]{}
	Bv=\dt(k*v) = \dt(k*\ell*f) = f\quad\mbox{in}\quad L^p(0,T;X),
$$
while also
$$
	B_nv = \dt(k*h_n*\ell*f) = \dt(k*\ell*h_n*f) = (h_n*f).
$$
This implies that \eqref{eq:2.5} holds.
In particular, putting $f=k$ with \eqref{eq:2.4}, we have
$$
	k_n \to k\quad\mbox{strongly in $L^1(0,T)$}\quad\mbox{as}\quad n \to \infty.
$$

Next we recall an important convexity inequality for the nonlocal time derivative operator.
The following inequality can be regarded as a nonlocal analogue of the classical chain rule $(H(u))' = H'(u)u'$.
\begin{lemma}
\label{Lemma:2.1}\cite{KSVZ}*{Corollary~6.1}
	Let $T >0$ and $I\subset\R$ be open.
	Let further $h \in W^{1,1}(0,T)$ be a nonnegative nonincreasing function, 
	$H \in C^1(I)$ be convex, $u_0 \in I$, and $u \in L^1(0,T)$ with $u(t) \in I$ for a.e. $t \in (0,T)$.
	Suppose that the functions $H(u), H'(u)u$, and $H'(u)(\dot{h}*u)$ belong to $L^1(0,T)$.
	Then, it holds that
	$$
		H'(u(t))\dt \left(h*[u(\cdot)-u_0]\right)(t) 
		\ge \dt\left(h*\bigl[H(u(\cdot))-H(u_0)\bigr]\right)(t)\quad \mbox{a.e.}\quad t \in (0,T).
	$$
\end{lemma}
Finally, applying the same arguments as in \cite{AN}*{Lemma 4.3} directly, 
we have the following.
\begin{lemma}
\label{Lemma:2.2}
	Let $T>0$, $p\in(1,\infty)$, and let $X$ be a Banach space.
	Assume that $u \in L^1(0,T; X)$ with $u(0)=u_0$ satisfies that
	$$
		k*[u-u_0] \in W^{1,p}(0,T; X)\,\,\,\mbox{with}\,\,\,\left(k*[u-u_0]\right)(0)=0,
	$$
	and
	$$
	\ell*\left\|\dt\left(k*[u-u_0]\right)\right\|_X^p \in L^\infty(0,T).
	$$
	Then $u \in C([0,T]; X)$ and $u(0)=u_0$ in $X$.
\end{lemma}
{\bf Proof.}
Put
$$
	w(t):=\dt(k*[u-u_0])(t)
$$
for $0<t<T$.
We first recall that the identity
\begin{equation}
\label{eq:2.6}
	u(t) = u_0 + (\ell*w)(t)
\end{equation}
holds for a.e. $t\in(0,T)$.
In fact, since $k*\ell=1$ and $(k*[u-u_0])(0)=0$, it holds that
$$
	\begin{aligned}
	u(t)-u_0
	& = \dt\left(1*[u-u_0]\right)(t)
	= \dt\left(\ell*k*[u-u_0]\right)(t) 
	\\
	& = (\ell*w)(t) + \ell(t)\left(k*[u-u_0]\right)(0) 
	= (\ell*w)(t)
	\end{aligned}
$$
for a.e. $t\in(0,T)$.

Next we prove $u \in C([0,T]; X)$ and $u(0)=u_0$ in $X$.
Let $0<s<t\le T$.
Then, it follows from \eqref{eq:2.6} that
\begin{equation}
\label{eq:2.7}
	\begin{aligned}
	\|u(t)-u(s)\|_X
	& = \left\|\left(\ell*w\right)(t)-\left(\ell*w\right)(s)\right\|_X
	\\
	& = \left\|\int_0^t\ell(t-\tau)w(\tau)\,\dee\tau 
	- \int_0^s\ell(s-\tau)w(\tau)\,\dee\tau\right\|_X
	\\
	& \le \int_0^s\left|\ell(t-\tau)-\ell(s-\tau)\right|\|w(\tau)\|_X\,\dee\tau
	+ \int_s^t\ell(t-\tau)\|w(\tau)\|_X\ \dee\tau
	\\
	& =: I(t,s) + J(t,s).
	\end{aligned}
\end{equation}
For the first term, by H\"{o}lder's inequality, we see that
$$
	I(t,s)
	\le \left(\int_0^s\left|\ell(t-\tau)-\ell(s-\tau)\right|\,\dee\tau\right)^{1-1/p}
	\left(\int_0^s\left|\ell(t-\tau)-\ell(s-\tau)\right|\|w(\tau)\|_X\,\dee\tau\right)^{1/p}.
$$
Here we recall that $\ell$ is nonnegative and nonincreasing, so it follows that
$$
	\left|\ell(t-\tau)-\ell(s-\tau)\right| = \ell(s-\tau)-\ell(t-\tau) \leq \ell(s-\tau)\quad {\rm for}\quad 0<s<t\le T.
$$
Thus, since $\ell\in L^1(0,T)$, we obtain
$$
	\begin{aligned}
	I(t,s)
	& \le \left(\int_0^s\left|\ell(t-\tau)-\ell(s-\tau)\right|\,\dee\tau\right)^{1-1/p}
	\left(\int_0^s\ell(s-\tau)\|w(\tau)\|_X\,\dee\tau\right)^{1/p} 
	\\
	& \le \left(\int_0^s\left|\ell(t-\tau)-\ell(s-\tau)\right|\,\dee\tau\right)^{1-1/p}
	\left\|\left(\ell*\|w\|_X^p\right)\right\|_{L^\infty(0,T)}^{1/p}
	\to 0\quad{\rm as}\quad s\to t.
	\end{aligned}
$$
On the other hand, for the second term in \eqref{eq:2.7}, we also get
$$
	\begin{aligned}
	J(t,s)
	& = \int_s^t\ell^{1-1/p}(t-\tau)\ell^{1/p}(t-\tau)\|w(\tau)\|_X\,\dee\tau 
	\\
	& \le \left(\int_s^t\ell(t-\tau)\,\dee\tau\right)^{1-1/p}
	\left(\int_0^t\ell(t-\tau)\|w(\tau)\|_X^p\,\dee\tau\right)^{1/p} 
	\\
	& \le \left(\int_s^t\ell(t-\tau)\,\dee\tau\right)^{1-1/p}
	\left\|\left(\ell*\|w\|_X^p\right)\right\|_{L^\infty(0,T)}^{1/p}
	\to 0\quad{\rm as}\quad s\to t.
	\end{aligned}
$$
These together with \eqref{eq:2.7} imply $u \in C((0,T];X)$.
In the same way, we deduce that
$$
	\begin{aligned}
	\|u(t)-u_0\|_X
	& = \left\|\left(\ell*\left[\dt\left(k*[u-u_0]\right)\right]\right)(t)\right\|_X 
	\\
	& \le \left(\ell*\left\|\dt\left(k*[u-u_0]\right)\right\|_X\right)(t)
	\\
	& \le \|\ell\|_{L^1(0,t)}^{1-1/p}
	\left\|\left(\ell*\left\|\dt\left(k*[u-u_0]\right)\right\|_X^p\right)\right\|_{L^\infty(0,T)}^{1/p}
	\to 0\quad{\rm as}\quad{t\to0^+}.
\end{aligned}
$$
Hence we see that $u(0)=u_0$ in $X$, and $u \in C([0,T];X)$.
$\Box$
%
	\section{$L^1$-theory for bounded domains}\label{section:3}
Consider the initial-boundary value problem
\begin{equation}
\label{eq:IBP}\tag{IBP}
	\begin{cases}
		\partial_t\left(k*[u-u_0]\right)=\Delta\Phi(u), & x \in \Omega,\quad t\in(0,T),\vspace{5pt}\\
		u=0, & x \in \partial\Omega,\quad t \in (0,T),\vspace{5pt}\\
		u(x,0)=u_0(x), & x \in \Omega,
	\end{cases}
\end{equation}
where $\Omega\subset\R^N$ is a bounded domain with smooth boundary and $T>0$ may be finite or infinite.
In this section, we show the existence and uniqueness of global weak solutions to \eqref{eq:IBP}
which satisfy improved energy estimates.
Furthermore, we develop the $L^1$-theory for a bounded domain.

\subsection{Improved energy estimates for weak solutions}
We first introduce the definition of a weak solution to problem \eqref{eq:IBP}.
\begin{definition}
\label{Definition:3.1}
	Let $T>0$.
	We say that $u \in L^2(0,T;L^2(\Omega))$ is a weak solution to problem \eqref{eq:IBP} in $\Omega\times[0,T)$ if
	\begin{itemize}
		\item[\rm (i)] $k*[u-u_0] \in H^1(0,T; H^{-1}(\Omega))$ with $\left(k*[u-u_0]\right)(0)=0$,
		\item[\rm (ii)] $\Phi(u) \in L^2(0,T; H^1(\Omega))$,
		\item[\rm (iii)] $u$ satisfies the identity
		\begin{equation}
		\label{eq:3.1}
			\begin{aligned}
			&
			\int_0^T\int_{\Omega}(k*[u(x,\cdot)-u_0(x)])(t)\partial_t\psi(x,t)\,\dee x\,\dee t 
			\\
			&\hspace{5cm}
			-\int_0^T\int_{\Omega}\nabla\Phi(u(x,t))\cdot\nabla\psi(x,t)\,\dee x\,\dee t = 0
			\end{aligned}
		\end{equation}
		for any $\psi \in H^1(0,T; H_0^1(\Omega))$ with $\psi\bigl|_{t=T}=0$.
	\end{itemize}
	Moreover, 
	we say that $u \in L_{\rm loc}^2(0,\infty;L^2(\Omega))$ is a global weak solution to problem \eqref{eq:IBP}
	in $\Omega\times[0,\infty)$ if above conditions {\rm (i)}--{\rm (iii)} hold for any fixed $T>0$.
\end{definition}
For the non-degenerate and non-singular cases, we recall the following results.
\begin{Proposition}
\label{Proposition:3.1}{\cites{VZ2, WWZ, Z4}}
	Let $u_0 \in L^2(\Omega)$. 
	Suppose that $\Phi\in C^1(\R)$ satisfies that there exist $c_1, c_2>0$ such that
	\begin{equation}
	\label{phi}\tag{H$\Phi$}
		c_1 \le \Phi'(r) \le c_2,\qquad r\in\R.
	\end{equation}
	Then problem \eqref{eq:IBP} possesses a weak solution $u \in L^2_{\rm loc}(0,\infty; H_0^1(\Omega))$.
	Furthermore, for any $T>0$, the solution $u$ satisfies
	$$
		\|k*[u-u_0]\|_{H^1(0,T; H^{-1}(\Omega))} + \|u\|_{L^2(0,T; H_0^1(\Omega))} 
		\le C\|u_0\|_{L^2(\Omega)},
	$$
	where $C$ depends only on $c_1$, $c_2$, $T$, and $\|\ell\|_{L^1(0,T)}$.
\end{Proposition}
\begin{proposition}
\label{Proposition:3.2}\cite{Z3}*{Theorem 4.2}
	Let $u_0 \in L^2(\Omega)$ and let $u$ be weak solutions to problem \eqref{eq:IBP}. 	
	Suppose that $\Phi \in C^1(\R)$ satisfies \eqref{phi}. 	
	Then it holds that
	$$
		\essinf_{x\in\Omega}(u_0) \le u(x,t) \le \esssup_{x\in\Omega}(u_0)\quad\mbox{a.e.}\quad x \in \Omega,\,\,\,t \in (0,T).
	$$
\end{proposition}
For the degenerate case, that is, the case $\Phi'(0)=0$ (e.g., $\Phi(r)=|r|^{m-1}r$, $m>1$),
Wittbold {\it et al.} \cite{WWZ}*{Theorem 6.1} proved the existence of bounded weak solutions to \eqref{eq:IBP}.
We next extend the above existence result to the singular case, that is,
the case where $\Phi'(r)\to\infty$ as $r\to0^+$ (e.g., $\Phi(r)=|r|^{m-1}r$, $m\in(0,1)$).
Furthermore, we prove the uniqueness of weak solutions not only for the degenerate case
but also for the singular case, and provide improved energy estimates.
\begin{theorem}
\label{Theorem:3.1}
	Assume $(K)$ and $(L)$.
	Let $u_0$ be a bounded function.
	Then the following assertions hold.
	\begin{itemize}
	\item[\rm(i)]
		Problem \eqref{eq:IBP} possesses a global weak solution 
		$u \in C([0,\infty);H^{-1}(\Omega))$ with $u(0)=u_0$ in $H^{-1}(\Omega)$ satisfying
		\begin{equation}
		\label{eq:3.2}
		\ell*\|\nabla\Phi(u)\|_{L^2(\Omega)}^2 \in L^\infty(0,\infty),\quad
		\ell*\|\partial_t(k*[u-u_0])\|_{H^{-1}(\Omega)}^2 \in L^\infty(0,\infty).
		\end{equation}
		Furthermore, for any $q\in[1,\infty]$, the solution $u$ also belongs to 
		$L^\infty(0,\infty\,;\,L^q(\Omega))$, and satisfies
		\begin{equation}
		\label{eq:3.3}
			\|u\|_{L^\infty(0,\infty\,;\,L^q(\Omega))} \le\|u_0\|_{L^q(\Omega)}.
		\end{equation}
	\item[\rm(ii)]
		Let $u, v$ be global weak solutions to problem \eqref{eq:IBP} with bounded initial data 
		$u_0,v_0$, respectively.
		Furthermore, let $\zeta$ be a nonnegative classical solution to problem
		\begin{equation}
		\label{eq:3.4}
			-\Delta\zeta = 1\quad{\rm in}\quad\Omega,\qquad
			\zeta=0\quad{\rm on}\quad\partial\Omega.
		\end{equation}
		Then it holds that
		\begin{equation}
		\label{eq:3.5}
			\|u(t)-v(t)\|_{L^1_\zeta(\Omega)} + \left(\ell*\|\Phi(u)-\Phi(v)\|_{L^1(\Omega)}\right)(t)
			\le \|u_0-v_0\|_{L^1_\zeta(\Omega)},\quad t\ge0,
		\end{equation}
		where $\|f\|_{L^1_\zeta(\Omega)}:=\|f\zeta\|_{L^1(\Omega)}$.
		In particular,
		\begin{equation}
		\label{eq:3.6}
			\int_{\Omega}[u(x,t)-v(x,t)]_+\zeta(x)\,\dee x 
			\le \int_{\Omega}[u_0(x)-v_0(x)]_+\zeta(x)\,\dee x,
			\quad t\ge0.
		\end{equation}
	\end{itemize}
\end{theorem}
\begin{remark}
\label{Remark:3.1}
	By the standard elliptic theory, it holds that $\zeta>0$ in $\Omega$ $($see, e.g., \cite[Chapter 6.6]{V2}$)$.
	By assertion~{\rm (ii)} of Theorem~$\ref{Theorem:3.1}$, since $\zeta>0$ in $\Omega$, we see that
	the solution $u$ to \eqref{eq:IBP} is unique, 
	and the comparison principle for weak solutions to \eqref{eq:IBP} holds,
	that is, if $u_0 \le v_0$ a.e. $x\in\Omega$, then $u \le v$ a.e. $x\in\Omega$ and $t>0$.
	In particular, if $u_0\ge0$, then $u\ge0$ a.e. $x\in\Omega$ and $t>0$.
\end{remark}
{\bf Proof of assertion~(i) of Theorem~\ref{Theorem:3.1}.}
Put $M:=\|u_0\|_{L^\infty(\Omega)}$+1.
Let $n\in\mathbb N$.
In order to avoid degeneracy or singularity, we first define
\begin{equation}
\label{eq:3.7}
	\tilde\Phi(r) := 
	\left\{
	\begin{array}{ll}
	\displaystyle{\Phi'(M)(r-M)+\Phi(M)}
	&\mbox{if}\quad r>M,
	\vspace{5pt}\\
	\displaystyle{\Phi(r)}
	&\mbox{if}\quad r\in[-M,M],
	\vspace{5pt}\\
	\displaystyle{\Phi'(M)(r+M)-\Phi(M)}
	&\mbox{if}\quad r<-M,
	\end{array}
	\right.
\end{equation}
and put
\begin{equation}
\label{eq:3.8}
	\tilde\Phi_n(r) := 
	\left\{
	\begin{array}{ll}
	\displaystyle{\tilde\Phi\left(r+\frac{1}{n}-\alpha_n\right)}
	&\mbox{if}\quad r>\alpha_n,
	\vspace{5pt}\\
	\displaystyle{\Phi'\left(\frac{1}{n}\right)r}
	&\mbox{if}\quad r\in[-\alpha_n,\alpha_n],
	\vspace{5pt}\\
	\displaystyle{\tilde\Phi\left(r-\frac{1}{n}+\alpha_n\right)}
	&\mbox{if}\quad r<-\alpha_n,
	\end{array}
	\right.
\end{equation}
where
$$
	\alpha_n:=\Phi\left(\frac{1}{n}\right)\left\{\Phi'\left(\frac{1}{n}\right)\right\}^{-1}.
$$
Then, by condition $(L)$ we see that $\tilde\Phi_n\in C^1(\R)$, $\alpha_n\to0$ as $n\to\infty$, and 
the function $\tilde\Phi_n(r)$ satisfies \eqref{phi}.
Consider the initial-boundary value problem
\begin{equation}
\label{IBP-n}\tag{${\rm IBP}_n$}
	\begin{cases}
		\partial_t\left(k*[u_n-u_{0}]\right) = \Delta \tilde\Phi_n(u_n), 
		& x\in\Omega,\,\,\,t>0,
		\vspace{5pt} \\
		u_n=0, 
		& x\in\partial\Omega,\,\,\,t>0,
		\vspace{5pt} \\
		u_n(x,0)=u_0(x), 
		& x\in\Omega.
	\end{cases}
\end{equation}
Since $u_0\in L^\infty(\Omega)$, namely $u_0\in L^2(\Omega)$, 
thanks to Proposition~\ref{Proposition:3.1}, 
there exists a weak solution $u_n \in L_{\rm loc}^2(0,\infty;H_0^1(\Omega))$.
Furthermore, it follows from Proposition~\ref{Proposition:3.2} that
\begin{equation}
\label{eq:3.9}
	\|u_n\|_{L^\infty(\Omega\times(0,\infty))} \le \|u_0\|_{L^\infty(\Omega)}.
\end{equation}

Next we construct a solution which belongs to $C([0,\infty);H^{-1}(\Omega))$.
By \cite{Z3}*{Lemma~3.1}, 
$u_n$ satisfies the time-regularized version of the weak formulation of problem \eqref{IBP-n},
that is, for any $\varphi \in H_0^1(\Omega)$,
\begin{equation}
\label{eq:3.10}
	\begin{aligned}
	&
	\int_{\Omega}\partial_t\left(k_j*[u_n(x,\cdot)-u_{0}(x)]\right)(t)\varphi(x)\,\dee x 
	\\
	&\qquad
	+\int_{\Omega}\left(h_j*\left[\nabla\tilde\Phi_n(u_n(x,\cdot))\right]\right)(t)\cdot\nabla\varphi(x)\,\dee x = 0,
	\quad {\rm a.e.}\,\,\,t>0,\quad j\in\N,
	\end{aligned}
\end{equation}
where $h_j$ is defined by \eqref{eq:2.3}, and $k_j = k*h_j$.
Taking in \eqref{eq:3.10} the test function $\varphi=\tilde\Phi_n(u_n)$,
we have
\begin{equation}
\label{eq:3.11}
	\begin{aligned}
	&
	\int_\Omega\partial_t(k_j*[u_n(x,\cdot)-u_0(x)])(t)\tilde\Phi_n(u_n(x,t))\,\dee x 
	\\
	&\qquad
	+ \int_{\Omega}\biggl(h_j*\biggl[\nabla\tilde\Phi_n(u_n(x,\cdot))\biggr]\biggr)(t)\cdot\nabla\tilde\Phi_n(u_n(x,t))\,\dee x=0,
	\quad {\rm a.e.}\,\,\,t>0.
	\end{aligned}
\end{equation}
We now define
$$
	\Psi_n(r) := \int_0^r\tilde \Phi_n(s)\,\dee s,\ \ r \in \R.
$$
Then $\Psi_n'=\tilde \Phi_n$ and $\Psi_n'' = \tilde \Phi_n' \ge 0$ hold, 
so $\Psi_n$ is convex.
Applying Lemma~\ref{Lemma:2.1} with $H(y)=\Psi_n(y)$, we have
$$
	\partial_t(k_j*[u_n(x,\cdot)-u_0(x)])(t)\tilde \Phi_n(u_n(x,t)) 
	\ge \partial_t\biggl(k_j*\biggl[\Psi_n(u_n(x,\cdot))-\Psi_n(u_0(x))\biggr]\biggr)(t)
$$
for a.e. $x \in \Omega$ and $t>0$.
This together with \eqref{eq:3.11} yields
\begin{equation}
\label{eq:3.12}
	\begin{aligned}
	&
	\int_\Omega\partial_t\biggl(k_j*\biggl[\Psi_n(u_n(x,\cdot))-\Psi_n(u_0(x))\biggr]\biggr)(t)\,\dee x
	\\
	&\qquad
	+ \int_\Omega\biggl(h_j*\biggl[\nabla\tilde\Phi_n(u_n(x,\cdot))\biggr]\biggr)(t)\cdot\nabla\tilde\Phi_n(u_n(x,t))\,\dee x
	\le 0,
	\quad{\rm a.e.}\,\,\,t>0.
	\end{aligned}
\end{equation}
Hence we convolve \eqref{eq:3.12} with the kernel $\ell$ to get
$$
	\begin{aligned}
	&
	\int_\Omega\biggl(h_j*\biggl[\Psi_n(u_n(x,\cdot))-\Psi_n(u_0(x))\biggr]\biggr)(t)\,\dee x 
	\\
	&\qquad
	+ \left(\ell*\left[\int_\Omega\biggl(h_j*\biggl[\nabla\tilde \Phi_n(u_n(x,\cdot))\biggr]\biggr)\cdot
	\nabla\tilde \Phi_n(u_n(x,\cdot))\,\dee x\right]\right)(t)
	\le 0,
	\quad {\rm a.e.}\,\,\,t>0.
	\end{aligned}
$$
Here we used facts that $k*\ell=1$,
$$
	\biggl(k_j*\biggl[\Psi_n(u_n)-\Psi_n(u_0)\biggr]\biggr)(0) = \left(k_j*[u_n-u_0]\right)(0) = 0,
$$
and
$$
	\ell*[\partial_t(k_j*[u_n-u_0])] = \partial_t(\ell*k*h_j*[u_n-u_0]) = (h_j*[u_n-u_0]).
$$
Since it follows from \eqref{eq:2.5} that 
$h_j*f \to f$ strongly in $L^1(0,T; L^1(\Omega))$ for any $f \in L^1(0,T;L^1(\Omega))$, 
we can take an appropriate subsequence $h_j*f\to f$ a.e. in $\Omega\times(0,T)$ if necessary,
and we obtain
\begin{equation}
\label{eq:3.13}
	\int_\Omega\biggl(\Psi_n(u_n(x,t))-\Psi_n(u_0(x))\biggr)\,\dee x
	+ \left(\ell*\left[\int_\Omega\left|\nabla\tilde \Phi_n(u_n(x,\cdot))\right|^2\,\dee x\right]\right)(t) 
	\le 0,
	\quad{\rm a.e.}\,\,\,t>0.
\end{equation}
This implies that
\begin{equation}
\label{eq:3.14}
	\left(\ell*\|\nabla\tilde \Phi_n(u_n)\|_{L^2(\Omega)}^2\right)(t) 
	\le \int_\Omega\Psi_n(u_0(x))\,\dee x,
	\,\,\,{\rm a.e.}\,\,t>0.
\end{equation}
On the other hand, choosing $\varphi=u_n$ as a test function in \eqref{eq:3.10} 
and applying Lemma~\ref{Lemma:2.1} with $H(y)=\frac{1}{2}y^2$, we have
$$
	\int_\Omega u_n^2(x,t)\,\dee x 
	+ 2\left(\ell*\left[\int_{\Omega}\tilde \Phi_n'(u_n(x,\cdot))|\nabla u_n(x,\cdot)|^2\,\dee x\right]\right)(t) 
	\le \int_\Omega u_0^2(x)\,\dee x,\quad{\rm a.e.}\,\,\,t>0.
$$
Since the second term is nonnegative, we drop this term to get
\begin{equation}
\label{eq:3.15}
	\int_\Omega u_n^2(x,t)\,\dee x 
	\le \int_\Omega u_0^2(x)\,\dee x,\quad{\rm a.e.}\,\,\, t>0.
\end{equation}
Let $w \in H_0^1(\Omega)$ satisfying $\|w\|_{H_0^1(\Omega)}=1$. 
Taking $\varphi=w$ as a test function in \eqref{eq:3.10}, we deduce that
$$
	\begin{aligned}
	&
	\int_\Omega\partial_t(k*[u_n(x,\cdot)-u_0(x)])(t)w(x)\,\dee x 
	\\
	& 
	= -\int_\Omega\nabla\tilde\Phi_n(u_n(x,t))\cdot\nabla w(x)\,\dee x 
	+ \int_\Omega\partial_t\biggl((k-k_j)*[u_n(x,\cdot)-u_0(x)]\biggr)(t)w(x)\,\dee x 
	\\
	& \hspace{2cm} 
	+ \int_\Omega\biggl\{\nabla\tilde{\Phi}_n(u_n(x,t)) 
	- \left(h_j*\biggl[\nabla\tilde \Phi_n(u_n(x,\cdot))\biggr]\right)(t)\biggr\}
	\cdot\nabla w(x)\,\dee x,
	\quad{\rm a.e.}\,\,\,t>0.
	\end{aligned}
$$
It follows from $k*[u_n-u_0](0)=0$ that
$$
	\partial_t(k_j*[u_n-u_0]) 
	= \partial_t(h_j*k*[u_n-u_0]) 
	= h_j*[\partial_t(k*[u_n-u_0])],
$$
so we see that
$$
	\begin{aligned}
	& 
	\int_\Omega\partial_t\left(k*[u_n(x,\cdot)-u_0(x)]\right)(t)w(x)\,\dee x 
	\\
	& 
	= -\int_\Omega\nabla\tilde\Phi_n(u_n(x,t))\cdot\nabla w(x)\,\dee x 
	+ \int_\Omega\biggl\{V(x,t) - \left(h_j*[V(x,\cdot)]\right)(t)\biggr\}\,\dee x,
	\,\,\,{\rm a.e.}\,\,t>0,
	\end{aligned}
$$
where
$$
	V = \partial_t\left(k*[u_n-u_0]\right)w 
	+ \nabla\tilde \Phi_n(u_n)\cdot\nabla w \in L_{\rm loc}^1(0,\infty;L^1(\Omega)).
$$
Then,
since it follows from \eqref{eq:2.5} that $h_j*V \to V$ a.e. $t>0$ as $j\to\infty$ 
by choosing an appropriate subsequence if necessary,
applying H\"{o}lder's inequality with $\|\nabla w\|_{L^2(\Omega)} \le \|w\|_{H_0^1(\Omega)} = 1$, 
we deduce that
\begin{equation}
\label{eq:3.16}
	\left|\int_\Omega\partial_t\left(k*[u_n(x,\cdot)-u_0(x)]\right)(t)w(x)\,\dee x\right|^2 
	\le \|\nabla\tilde \Phi_n(u_n(t))\|_{L^2(\Omega)}^2,
	\,\,\,{\rm a.e.}\,\,t>0.
\end{equation}
Hence, by \eqref{eq:3.14} and 
$$
	\left\|\partial_t\left(k*[u_n-u_0]\right)(\cdot)\right\|_{H^{-1}(\Omega)} 
	:= \sup_{\|w\|_{H_0^1(\Omega)}=1}
	\left|\int_{\Omega}\partial_t\left(k*[u_n(x,\cdot)-u_0(x)]\right)(t)w(x)\,\dee x\right|,
$$
we convolve \eqref{eq:3.16} with the kernel $\ell$ to obtain 
\begin{equation}
\label{eq:3.17}
	\begin{aligned}
	&
	\left(\ell*\left\|\partial_t\left(k*[u_n-u_0]\right)(\cdot)\right\|_{H^{-1}(\Omega)}^2\right)(t)
	\\
	&
	\le \left(\ell*\|\nabla\tilde\Phi_n(u_n(\cdot))\|_{L^2(\Omega)}^2\right)(t) 
	\le \|\Psi_n(u_0)\|_{L^1(\Omega)},
	\,\,\,{\rm a.e.}\,\,t>0.
	\end{aligned}
\end{equation}
Convolving with the kernel $k$, for any $T>0$, we have
$$
	\begin{aligned}
	&
	\int_0^T\left\|\partial_t\left(k*[u_n-u_0]\right)(t)\right\|_{H^{-1}(\Omega)}^2\,\dee t 
	\\
	&
	\le \int_0^T\|\nabla\tilde \Phi_n(u_n(t))\|_{L^2(\Omega)}^2\,\dee t 
	\le \|k\|_{L^1(0,T)}\|\Psi_n(u_0)\|_{L^1(\Omega)}.
	\end{aligned}
$$
Therefore, by \eqref{eq:3.14}, \eqref{eq:3.15}, \eqref{eq:3.17}, 
and using the weak compactness argument, 
we obtain $u \in L^\infty(0,\infty;L^2(\Omega))$ satisfying
\begin{equation}
\label{eq:3.18}
	\ell*\|\nabla v_1\|_{L^2(\Omega)}^2 \in L^\infty(0,\infty),
	\quad
	\ell*\left\|\partial_tv_2\right\|_{H^{-1}(\Omega)}^2 \in L^\infty(0,\infty),
\end{equation}
and a subsequence $\{u_{n_m}\}$ such that
$$
	\begin{aligned}
	u_{n_m}& \to u\quad\mbox{weakly in $L_{\rm loc}^2(0,\infty;L^2(\Omega))$ as $m\to\infty$},
	\\
	\tilde\Phi_{n_m}(u_{n_m})& \to v_1\quad
	\mbox{weakly in $L_{\rm loc}^2(0,\infty;H_0^1(\Omega))$ as $m\to\infty$},
	\\
	k*[u_{n_m}-u_0]& \to v_2\quad
	\mbox{weakly in $H_{\rm loc}^1(0,\infty;H^{-1}(\Omega))$ as $m\to\infty$}.
	\end{aligned}
$$
Here we see that $v_2=(k*[u-u_0])$ in $H_{\rm loc}^1(0,\infty;H^{-1}(\Omega))$.
Indeed, by \cite{WWZ}*{Theorem~3.3} with $p=2$, $V=L^2(\Omega)$, and $H=H^{-1}(\Omega)$, 
it holds that there exists a subsequence $\{u_{n_m}\}$ such that
$$
	u_{n_m} \to u\quad\mbox{strongly in $L_{\rm loc}^2(0,\infty;H^{-1}(\Omega))$ as $m\to\infty$}.
$$
Therefore, by Young's inequality, for any fixed $T>0$, we obtain
$$
	\begin{aligned}
	& 
	\left|\int_0^T\int_\Omega(k*[u_{n_m}(x,\cdot)-u_0(x)])(t)\,\partial_t\psi(x,t)\,\dee x\,\dee t\right.
	\\
	&\hspace{3cm}
	\left.-\int_0^T\int_\Omega(k*[u(x,\cdot)-u_0(x)])(t)\,\partial_t\psi(x,t)\,\dee x\,\dee t\right| 
	\\
	& 
	\le \|k*[u_{n_m}-u]\|_{L^2(0,T;H^{-1}(\Omega))}\|\partial_t\psi\|_{L^2(0,T;H_0^1(\Omega))} 
	\\
	&
	\le \|k\|_{L^1(0,T)}\|u_{n_m}-u\|_{L^2(0,T;H^{-1}(\Omega))}
	\|\partial_t\psi\|_{L^2(0,T;H_0^1(\Omega))}
	\to 0\,\,\,{\rm as}\,\, m\to\infty
\end{aligned}
$$
for all $\psi\in H^1(0,T;H_0^1(\Omega))$ with $\psi\bigl|_{t=T}=0$.
On the other hand, applying the same argument as in the proof of \cite{WWZ}*{Theorem~6.1}
with \eqref{eq:3.7}, \eqref{eq:3.8}, and \eqref{eq:3.9},
we also see that $v_1=\Phi(u)$ in the sense of distributions,
namely, $\tilde \Phi_{n_m}(u_{n_m}) \to \Phi(u)$ weakly in $L_{\rm loc}^2(0,\infty;H_0^1(\Omega))$ 
as $m\to\infty$. (See also, e.g., \cite{S}*{Section~II~2}).
Thus, it follows from \eqref{eq:3.17} that 
$u$ satisfies the weak formulation \eqref{eq:3.1} and the estimate \eqref{eq:3.2}.
From \eqref{eq:3.18} and Lemma \ref{Lemma:2.2} with $X=H^{-1}(\Omega)$ and $p=2$, 
we obtain $u \in C([0,\infty);H^{-1}(\Omega))$ with $u(0)=u_0$ in $H^{-1}(\Omega)$.

Finally we prove \eqref{eq:3.3}.
For the case of  $p=\infty$, 
\eqref{eq:3.3} follows from \eqref{eq:3.9}.
Hence, we focus on the cases $1\le p<\infty$.
Let $p\in[1,\infty)$ and $u$ be a weak solution to \eqref{eq:IBP}.
We define the $C^\infty$-function
\begin{equation}
\label{eq:3.19}
	H_{\E}(y) := (y^2 + {\E}^2)^{\frac{p}{2}} - \E^p,\,\,\, \E>0,\,\, y\in\R.
\end{equation}
Then, it holds that, for $y\in\R$,
$$
	H_{\E}'(y) = py(y^2 + {\E}^2)^{\frac{p-2}{2}},\quad
	H_{\E}''(y) = \{(p-1)y^2+{\E}^2\}(y^2 + {\E}^2)^{\frac{p-4}{2}},
$$
and
\begin{equation}
\label{eq:3.20}
	|H_{\E}(y)| \le 2^{p-1}(|y|^p+{\E}^p),\qquad 
	H_{\E}(y) \to |y|^p\quad \mbox{as}\quad \E\to0^+.
\end{equation}
Thus we see that $H_{\E}$ is convex, as well as that $H_{\E}'$ and $H_{\E}''$ are bounded for $|y|\le M$.
Note that $H_{\E}'(u_n(t)) \in H_0^1(\Omega) \cap L^\infty(\Omega)$ for a.e. $t>0$,
and using this function as a test function in the time-regularized weak formulation \eqref{eq:3.10}, we get
$$
	\begin{aligned}
	&
	\int_\Omega\partial_t(k_j*[u_n(x,\cdot)-u_{0}(x)])(t)H_{\E}'(u_n(x,t))\,\dee x 
	\\
	&\qquad
	+ \int_\Omega\left(h_j*\left[\nabla\tilde{\Phi}_n(u_n(x,\cdot))\right]\right)(t)\cdot\nabla u_n(x,t)H_{\E}''(u_n(x,t))\,\dee x = 0,
	\quad{\rm a.e.}\,\,\, t>0.
	\end{aligned}
$$
By Lemma~\ref{Lemma:2.1} with $H(y)=H_{\E}(y)$ we deduce that
$$
	\partial_t(k_j*[u_n(x,\cdot)-u_{0}(x)])(t)H_{\E}'(u(x,t))
	\ge \partial_t\biggl(k_j*\biggl[H_{\E}(u_n(x,\cdot))-H_{\E}(u_{0}(x))\biggr]\biggr)(t)
$$
for a.e. $x\in\Omega$ and $t>0$.
Applying the same argument as in the proof of \eqref{eq:3.13}, we see that
$$
	\begin{aligned}
	&
	\int_\Omega H_{\E}(u_n(x,t))\,\dee x 
	+\left(\ell*\left[\int_{\Omega}\tilde{\Phi}_n'(u_n(x,\cdot))|\nabla u_n(x,\cdot)|^2
	H_{\E}''(u_n(x,\cdot))\,\dee x\right]\right)(t)
	\\
	&\qquad\qquad 
	\le \int_{\Omega}H_{\E}(u_0(x))\,\dee x,
	\quad{\rm a.e.}\,\,\, t>0.
	\end{aligned}
$$
Since the second term in the left-hand side is nonnegative,
it holds that
$$
	\int_\Omega H_{\E}(u_n(x,t))\,\dee x \le \int_{\Omega}H_{\E}(u_0(x))\,\dee x,
	\,\,\,{\rm a.e.}\,\, t>0.
$$
By \eqref{eq:3.20}, we can apply the dominated convergence theorem to have
$$
	\int_\Omega H_{\E}(u_0(x))\,\dee x \to \|u_0\|_{L^p(\Omega)}^p,
	\qquad
	\int_\Omega H_{\E}(u_n(x,t))\,\dee x \to \|u_n(t)\|_{L^p(\Omega)}^p,\quad {\rm a.e.}\,\,\,t>0,
$$
as $\E\to0^+$ for all $n\in\N$.
This implies that \eqref{eq:3.3} for the case $p\in[1,\infty)$.
Thus assertion~(i) of Theorem~\ref{Theorem:3.1} follows.
$\Box$
\medskip

\noindent
{\bf Proof of assertion~(ii) of Theorem~\ref{Theorem:3.1}.}
Let $u, v$ be global weak solutions to problem \eqref{eq:IBP} with bounded initial data 
$u_0,v_0$, respectively,
and let $\zeta$ be a nonnegative classical solution to problem \eqref{eq:3.4}.
Set
$$
	U=u-v,\quad
	U_0=u_0-v_0,\quad
	V=\Phi(u)-\Phi(v).
$$
We also denote $H_{\E}$ defined as \eqref{eq:3.19} with $p=1$.
Taking $\varphi=H_{\E}'(V(t))\zeta \in H_0^1(\Omega)$ 
as a test function in the weak formulations \eqref{eq:3.10} for $u$, $v$, 
and taking the difference of the identities, we obtain
\begin{equation}
\label{eq:3.21}
	\begin{aligned}
	&
	\int_{\Omega}\partial_t(k_{j}*[U(x,\cdot)-U_0(x)])(t)H_{\E}'(V(x,t))\zeta(x)\,\dee x 
	\\
	&\qquad
	+ \int_{\Omega}(h_{j}*[\nabla V(x,\cdot)])(t)\cdot\nabla V(x,t)H_{\E}''(V(x,t))\zeta(x)\, \dee x 
	\\
	&\qquad\qquad
	+ \int_{\Omega}(h_{j}*[\nabla V(x,\cdot)])(t)\cdot\nabla\zeta(x)H_{\E}'(V(x,t))\,\dee x = 0,
	\quad{\rm a.e.}\,\,\, t>0.
	\end{aligned}
\end{equation}
Furthermore, for the first term in the left-hand side, we have
$$
	\begin{aligned}
	& \int_\Omega\partial_t(k_j*[U(x,\cdot)-U_0(x)])(t)H_{\E}'(V(x,t))\zeta(x)\,\dee x 
	\\
	&\qquad
	= \int_\Omega\partial_t(k_j*[U(x,\cdot)-U_0(x)])(t)H_{\E}'(U(x,t))\zeta(x)\,\dee x 
	\\
	&\qquad\qquad  
	+ \int_\Omega\partial_t(k_j*[U(x,\cdot)-U_0(x)])(t)\biggl\{H_{\E}'(V(x,t)) - H_{\E}'(U(x,t))\biggr\}\zeta(x)\,\dee x 
	\\
	&\qquad
	=: \int_\Omega\partial_t(k_j*[U(x,\cdot)-U_0(x)])(t)H_{\E}'(U(x,t))\zeta(x)\,\dee x + I_{j, \E}(t),
	\quad{\rm a.e.}\,\,\, t>0.
	\end{aligned}
$$
Applying Lemma \ref{Lemma:2.1} with $H(y)=H_{\E}(y)$, we deduce that
$$
	\partial_t(k_j*[U(x,\cdot)-U_0(x)])(t)H_{\E}'(U(x,t))\zeta(x)
	\ge \partial_t\biggl(k_j*\biggl[H_{\E}(U(x,\cdot))-H_{\E}(U_0(x))\biggr]\biggr)(t)\zeta(x)
$$
for a.e. $x\in\Omega$ and $t>0$.
Then, by convolving \eqref{eq:3.21} with the kernel $\ell$, we see that
$$
	\begin{aligned}
	&
	\int_\Omega\biggl(h_j*\biggl[H_{\E}(U(x,\cdot))-H_{\E}(U_0(x))\biggr]\biggr)(t)\zeta(x)\,\dee x 
	+ \left(\ell*I_{j,\E}\right)(t) 
	\\
	&\qquad
	+ \left(\ell*\left[\int_\Omega(h_j*[\nabla V(x,\cdot)])(t)\cdot\nabla V(x,t)H_{\E}''(V(x,\cdot))\zeta(x)\,\dee x\right]\right)(t) 
	\\
	&\qquad\qquad
	+ \left(\ell*\left[\int_\Omega(h_j*[\nabla V(x,\cdot)])(t)\cdot\nabla\zeta(x)H_{\E}'(V(x,\cdot))\,\dee x\right]\right)(t) \le 0,
	\quad{\rm a.e.}\,\,\, t>0.
	\end{aligned}
$$
Here we can take an appropriate subsequence $h_j*f\to f$ a.e. in $\Omega\times(0,T)$ as $j\to\infty$
if necessary, and we obtain
\begin{equation}
\label{eq:11}
	\begin{aligned}
	&
	\int_\Omega\biggl(H_{\E}(U(x,t))-H_{\E}(U_0(x))\biggr)\zeta(x)\,\dee x 
	+ \left(\ell*\left[\int_\Omega|\nabla V(x,\cdot)|^2H_{\E}''(V(x,\cdot))\zeta(x)\,\dee x\right]\right)(t) 
	\\
	&\hspace{1cm}
	+ \left(\ell*\left[\int_{\Omega}\nabla H_{\E}(V(x,\cdot))\cdot\nabla\zeta(x)\,\dee x\right]\right)(t)
	 \le \left|\left(\ell*I_{\E}\right)(t)\right|,
	 \quad{\rm a.e.}\,\,\,t>0,
	\end{aligned}
\end{equation}
where
$$
	\begin{aligned}
	I_{\E}(t)
	&
	:= \lim_{j\to\infty}I_{j,\E}(t)
	\\
	&
	=\int_\Omega\partial_t(k*[U(x,\cdot)-U_0])(t)\bigg\{H_{\E}'(V(x,t)) - H_{\E}'(U(x,t))\bigg\}\zeta(x)\,\dee x,
	\quad{\rm a.e.}\,\,\,t>0.
	\end{aligned}
$$
Since $\zeta>0$ in $\Omega$, the second term in the left-hand side is nonnegative.
Furthermore, for the third term in the left-hand side, by \eqref{eq:3.4}, we deduce that
$$
	\int_\Omega\nabla H_{\E}(V(x,t))\cdot\nabla\zeta(x)\,\dee x
	= -\int_\Omega H_{\E}(V(x,t))\Delta\zeta(x)\,\dee x 
	= \int_{\Omega}H_{\E}(V(x,t))\,\dee x,
	\quad{\rm a.e.}\,\,\,t>0.
$$
On the other hand, it holds that $\ell*I_{\E} \to 0$ as $\E\to0^+$.
In fact, since $|H_{\E}'(y)| \le 1$ for all $y\in\R$, 
we have an estimate, which is uniform with respect to $\E>0$,
$$
	\begin{aligned}
	&
	|(\ell*I_{\E})(t)|
	\\
	& 
	\le \left(\ell*\left[\int_\Omega|H_{\E}'(V(x,\cdot)) - H_{\E}'(U(x,\cdot))|
	|\partial_t(k*[U(x,\cdot)-U_0(x)])(s)||\zeta(x)|\,\dee x\right]\right)(t) 
	\\
	& 
	\le 2\|\ell\|_{L^1(0,t)}^{\frac{1}{2}}
	\left\|\ell*\|\partial_t(k*[U-U_0])\|_{H^{-1}(\Omega)}^2\right\|_{L^\infty(0,t)}^{\frac{1}{2}}
	\|\zeta\|_{H_0^1(\Omega)} < \infty,
	\quad t>0.
	\end{aligned}
$$
We note that $H_{\E}'(y) \to {\rm sign}(y)$ as $\E\to0^+$ for $y\in\R$. 
Besides, by the condition ($L$), we see that
$$
	{\rm sign}(U) = {\rm sign}(u-v) = {\rm sign}\left(\Phi(u)-\Phi(v)\right) = {\rm sign}(V).
$$
Then it follows by the dominated convergence theorem that
\begin{equation}
\label{eq:12}
	\ell*I_{\E} 
	\to 
	\ell*\left[\int_\Omega\partial_t(k*[u-u_0])\bigg({\rm sign}(V) - {\rm sign}(U)\bigg)\zeta\,\dee x\right] = 0
\end{equation}
as $\E\to0^+$.
Since $H_{\E}(y) \to |y|$ as $\E\to0^+$ for all $y\in\R$, 
applying the dominated convergence theorem again, we see that
$$
	\int_\Omega H_{\E}(U(x,t))\zeta(x)\,\dee x 
	\to \|U(t)\|_{L^1_\zeta(\Omega)} = \|u(t)-v(t)\|_{L^1_\zeta(\Omega)},
	\quad{\rm a.e.}\,\,\, t>0,
$$
and
$$
	\int_\Omega H_{\E}(U_0)\zeta(x)\,\dee x 
	\to \|u_0-v_0\|_{L^1_\zeta(\Omega)},
$$
as $\E\to0^+$.
Hence we obtain
$$
	\|u(t)-v(t)\|_{L^1_\zeta(\Omega)} 
	+ \left(\ell*\|\Phi(u)-\Phi(v)\|_{L^1(\Omega)}\right)(t) 
	\le \|u_0-v_0\|_{L^1_\zeta(\Omega)},
	\quad{\rm a.e.}\,\,\, t>0.
$$
We now recall that $u \in C([0,\infty); H^{-1}(\Omega))$ with $u(0)=u_0$ in $H^{-1}(\Omega)$, 
so we observe that $u\zeta \in C([0,\infty);L^1(\Omega))$ 
with $u(0)\zeta = u_0\zeta$ in $L^1(\Omega)$.
Furthermore, since it follows from $(L)$ and \eqref{eq:3.3} that
$\|\Phi(u)\|_{L^\infty(\Omega\times(0,\infty))} \le \Phi(\|u_0\|_{L^\infty(\Omega)})$,
by $(K)$ we see that
$\ell*\|\Phi(u)\|_{L^1(\Omega)} \in C([0,\infty))$ 
with $\left(\ell*\|\Phi(u)\|_{L^1(\Omega)}\right)(0)=0$.
Thus we have \eqref{eq:3.5}.
Moreover, put
$$
	\tilde{H}_{\E}(y):=\sqrt{[y]_+^2+{\E}^2} -\E,\ \ \E>0,\ \ y\in\R.
$$
Then since it follows that $\tilde{H}_{\E}(y) \to [y]_+$ as $\E\to0^+$ for $y\in\R$,
applying the same argument as in the proof of \eqref{eq:3.5} 
with $\tilde{H}_{\E}$ instead of $H_{\E}$, 
we obtain \eqref{eq:3.6}.
Thus assertion~(ii) of Theorem~\ref{Theorem:3.1} follows,
and the proof of Theorem~\ref{Theorem:3.1} is complete.
$\Box$
\bigskip

\subsection{Existence and uniqueness of $L^1$-solutions for bounded domain case}
In this subsection, we present the $L^1$-theory for problem \eqref{eq:IBP}.
We first define an $L^1$-solution to problem \eqref{eq:IBP}.
\begin{definition}
\label{Definition:3.2}
	Let $u_0\in L^1(\Omega)$.
	We say that $u\in L^\infty(0,\infty; L^1(\Omega))$ is an $L^1$-solution to problem \eqref{eq:IBP} if
	\begin{itemize}
		\item[\rm (i)] 
		$\Phi(u) \in L_{\rm loc}^1(0,\infty; L^1(\Omega))$ 
		with $\ell*\|\Phi(u)\|_{L^1(\Omega)} \in L^\infty(0,\infty)$,
		\item[\rm (ii)] 
		$u$ satisfies the identity
		\begin{equation}
		\label{eq:3.22}
			\int_\Omega(u(x,t)-u_0(x))\phi(x)\,\dee x 
			= \left(\ell*\left[\int_{\Omega}\Phi(u(x,\cdot))\Delta\phi(x)\,\dee x\right]\right)(t),
			\quad{\rm a.e.}\,\,\,t > 0,
		\end{equation}
		for any $\phi\in C_0^{\infty}(\Omega)$.
	\end{itemize}
\end{definition}
Under this definition,
we obtain the following existence and uniqueness result for $L^1$-solutions to problem \eqref{eq:IBP}
with initial data in $L^1$, which are not necessarily bounded.
\begin{theorem}
\label{Theorem:3.2}
	Assume that $(K)$ and $(L)$ hold.
	Let $u_0 \in L^1(\Omega)$.
	Then problem \eqref{eq:IBP} possesses a unique $L^1$-solution $u$
	satisfying
		$$
			\|u\|_{L^\infty(0,\infty; L^1(\Omega))} \leq \|u_0\|_{L^1(\Omega)}.
		$$
	Furthermore, it holds that
		$$
		\begin{aligned}
		u\in C([0,\infty); L^1_\zeta(\Omega))\quad
		&{\rm with}\quad
		u(0) = u_0\,\,\, {\rm in}\,\,\, L^1_\zeta(\Omega),
		\\
		\ell*\|\Phi(u)\|_{L^1(\Omega)} \in C([0,\infty))\quad
		&{\rm with}\quad
		 (\ell*\|\Phi(u)\|_{L^1(\Omega)})(0)=0,
		 \end{aligned}
		 $$
	where $\zeta$ is the nonnegative classical solution to problem \eqref{eq:3.4}.
	Moreover,	
	let $u, v$ be $L^1$-solutions to problem \eqref{eq:IBP} with initial data 
	$u_0,v_0\in L^1(\Omega)$, respectively.
	Then \eqref{eq:3.5} and \eqref{eq:3.6} hold.
\end{theorem}
\begin{remark}
\label{Remark:3.2}
	By Theorem~$\ref{Theorem:3.2}$ we see that
	the comparison principle for $L^1$-solutions to \eqref{eq:IBP} holds,
	that is, if $u_0 \le v_0$ a.e. $x\in\Omega$, then $u \le v$ a.e. $x\in\Omega$ and $t>0$.
	In particular, if $u_0\ge0$, then $u\ge0$ a.e. $x\in\Omega$ and $t>0$.
\end{remark}
{\bf Proof.}
We first show the existence of $L^1$-solutions to problem \eqref{eq:IBP}.
Let $n,m \in \N$, and
$$
	u_0^{n,m} =
	\left\{
	\begin{array}{lll}
		n
		& {\rm if} \quad u_0>n,
		\vspace{5pt}\\
		u_0
		& {\rm if} \quad u_0 \in [-m,n],
		\vspace{5pt}\\
		-m
		& {\rm if} \quad u_0<-m.
	\end{array}
	\right.
$$
Consider the approximate problem
\begin{equation}\label{eq:3.23}
	\begin{cases}
		\partial_t\left(k*[u_{n,m}-u_{0}^{n,m}]\right)=\Delta\Phi(u_{n,m}), & x\in\Omega,\ \ t>0,
		\vspace{5pt}\\
		u_{n,m}=0, & x \in \partial\Omega,\ \ t>0,
		\vspace{5pt}\\
		u_{n,m}(x,0)=u_{0}^{n,m}(x), & x\in\Omega.
	\end{cases}
\end{equation}
By Theorem~\ref{Theorem:3.1}, we see that 
problem~\eqref{eq:3.23} possesses a unique global weak solution
\begin{equation}
\label{eq:3.24}
	u_{n,m} \in C([0,\infty);H^{-1}(\Omega)) \cap L^\infty(0,\infty;L^2(\Omega))
	\,\,\,{\rm with}\,\,\, u_{n,m}(0)=u_{0}^{n,m}\,\,\,{\rm in}\,\,\, H^{-1}(\Omega)
\end{equation}
satisfying
$$
	-m \le u_{n,m} \le n,\quad{\rm a.e.}\,\,\,x\in\Omega,\,\,t>0,
$$
and
$$
	\ell*\|\nabla\Phi(u_{n,m})\|_{L^2(\Omega)}^2 \in L^\infty(0,\infty),
	\qquad
	\ell*\left\|\partial_t\left(k*[u_{n,m}-u_{0}^{n,m}]\right)\right\|_{H^{-1}(\Omega)}^2 \in L^\infty(0,\infty).
$$
Furthermore, for each fixed $n\in\N$, since $u_{0}^{n,m} \ge u_{0}^{n,m+1}$ a.e. in $\Omega$,
it follows from Remark~\ref{Remark:3.1} that
\begin{equation}
\label{eq:3.25}
	u_{n,m}(x,t) \ge u_{n,m+1}(x,t),\quad{\rm a.e.}\,\,\, x\in\Omega,\,\,t>0,
\end{equation}
for all $m\in\N$.
Moreover, by \eqref{eq:3.3} with $q=1$ and the definition of $u_{0}^{n,m}$ we deduce that
\begin{equation}
\label{eq:3.26}
	\|u_{n,m}\|_{L^\infty(0,\infty;L^1(\Omega))} \le \|u_{0}^{n,m}\|_{L^1(\Omega)} \le \|u_0\|_{L^1(\Omega)}.
\end{equation}
Thus we see that $\{u_{n,m}\}_m$ is monotone decreasing and uniformly bounded in $L^\infty(0,\infty;L^1(\Omega))$.

On the other hand, let $m^*\in\N$, it follows from \eqref{eq:3.5} that
\begin{equation}
\label{eq:3.27}
	\begin{aligned}
	\|u_{n,m}(t)-u_{n,m^*}(t)\|_{L^1_\zeta(\Omega)} 
	+\left(\ell*\|\Phi(u_{n,m})-\Phi(u_{n,m^*})\|_{L^1(\Omega)}\right)(t) \\
	\le \|u_{0}^{n,m}-u_{0}^{n,m^*}\|_{L^1_\zeta(\Omega)},\quad t\ge0.
	\end{aligned}
\end{equation}
This implies that
\begin{equation}
\label{eq:3.28}
	\begin{aligned}
	&
	\sup_{t\ge0}\,\|u_{n,m}(t)-u_{n,m^*}(t)\|_{L^1_\zeta(\Omega)} \to 0\quad{\rm as}\,\,\, m,m^*\to\infty,
	\\
	&
	\sup_{t\ge0}\,\left|\left(\ell*\|\Phi(u_{n,m})\|_{L^1(\Omega)}\right)(t)
	-\left(\ell*\|\Phi(u_{n,m^*})\|_{L^1(\Omega)}\right)(t)\right| \to 0\quad{\rm as}\,\,\,m,m^*\to\infty.
	\end{aligned}
\end{equation}
Furthermore, for any fixed $T>0$, convolving \eqref{eq:3.27} with the kernel $k$ and using $k*\ell=1$, 
we have
$$
	\|\Phi(u_{n,m})-\Phi(u_{n,m^*})\|_{L^1(0,T;L^1(\Omega))} 
	\le \|k\|_{L^1(0,T)}\|u_{0}^{n,m}-u_{0}^{n,m^*}\|_{L^1_\zeta(\Omega)},
$$
and it holds that
\begin{equation}
\label{eq:3.29}
	\left\|\Phi(u_{n,m})-\Phi(u_{n,m^*})\right\|_{L^1(0,T;L^1(\Omega))} \to 0\quad{\rm as}\,\,\,m,m^*\to\infty.
\end{equation}
Then, by \eqref{eq:3.24}, \eqref{eq:3.25}, \eqref{eq:3.26}, \eqref{eq:3.28}, and \eqref{eq:3.29}, 
we see that there exists a
$$
	u_n \in L^\infty(0,\infty;L^1(\Omega))
$$
with $u_n(0)=u_0^n$ in $L^1_\zeta(\Omega)$
such that
$$
	\begin{aligned}
	u_{n,m} \to u_n\quad &\mbox{strongly in $L^\infty(0,\infty;L^1(\Omega))$ as $m\to\infty$},
	\\
	u_{n,m} \to u_n\quad &\mbox{strongly in $C([0,\infty);L^1_\zeta(\Omega))$ as $m\to\infty$},
	\\
	\Phi(u_{n,m}) \to w_1^n\quad &\mbox{strongly in $L_{\rm loc}^1(0,\infty;L^1(\Omega))$ as $m\to\infty$},
	\\
	\ell*\|\Phi(u_{n,m})\|_{L^1(\Omega)} \to w_2^n\quad &\mbox{strongly in $C([0,\infty))$ as $m\to\infty$},
	\end{aligned}
$$
where $u_0^n=\min\{u_0,n\}$.
Applying the standard theory for nonlinear diffusion equations 
(see, e.g., \cite{V2}*{Section~5}) with condition $(L)$, 
we see that
$w_1^n=\Phi(u_n)$ and $w_2^n=\ell*\|\Phi(u_n)\|_{L^1(\Omega)}$, and that
the function $u_n$ satisfies
\begin{equation}
\label{eq:3.32}
	\|u_n\|_{L^\infty(0,\infty;L^1(\Omega))} \le \|u_0^n\|_{L^1(\Omega)} \le \|u_0\|_{L^1(\Omega)},
\end{equation}
and estimates \eqref{eq:3.5} and \eqref{eq:3.6}.
Since $u_0^n \le u_0^{n+1}$ a.e. in $\Omega$, we have
\begin{equation*}
	u_n(x,t) \le u_{n+1}(x,t), \quad{\rm a.e.}\,\,\, x\in\Omega,\,\,t>0,
\end{equation*}
for all $n\in\N$, which together with \eqref{eq:3.32} implies that
$\{u_n\}_n$ is monotone increasing and uniformly bounded in $L^\infty(0,\infty;L^1(\Omega))$.
By repeating the same argument, we see that there exists a $u \in L^\infty(0,\infty;L^1(\Omega))$ with $u(0)=u_0$ in $L_{\zeta}^1(\Omega)$ such that
$$
	\begin{aligned}
	u_{n} \to u\quad &\mbox{strongly in $L^\infty(0,\infty;L^1(\Omega))$ as $n\to\infty$},
	\\
	u_{n} \to u\quad &\mbox{strongly in $C([0,\infty);L^1_\zeta(\Omega))$ as $n\to\infty$},
	\\
	\Phi(u_n) \to \Phi(u)\quad &\mbox{strongly in $L_{\rm loc}^1(0,\infty;L^1(\Omega))$ as $n\to\infty$},
	\\
	\ell*\|\Phi(u_{n})\|_{L^1(\Omega)} \to \ell*\|\Phi(u)\|_{L^1(\Omega)}\quad &\mbox{strongly in $C([0,\infty))$ as $n\to\infty$}.
	\end{aligned}
$$
Furthermore, it also holds that $(\ell*\|\Phi(u)\|_{L^1(\Omega)})(0)=0$.
Indeed, since
$$
	\|\Phi(u_{n,m})\|_{L^\infty(\Omega\times(0,\infty))} \leq \Phi\bigl(\|u_0^{n,m}\|_{L^\infty(\Omega)}\bigr), \quad \mbox{a.e. $x\in\Omega$, $t>0$},
$$
it follows that
$$
	\begin{aligned}
	\left(\ell*\|\Phi(u)\|_{L^1(\Omega)}\right)(t)
	& 
	\le \left(\ell*\|\Phi(u)-\Phi(u_{n,m})\|_{L^1(\Omega)}\right)(t) + \left(\ell*\|\Phi(u_{n,m})\|_{L^1(\Omega)}\right)(t) 
	\\
	& 
	\le \|u_0-u_{0}^{n,m}\|_{L^1_\zeta(\Omega)} + \|\ell\|_{L^1(0,t)}\Phi\bigl(\|u_0^{n,m}\|_{L^\infty(\Omega)}\bigr)|\Omega|
	\to 0
	\end{aligned}
$$
as $t\to0^+$ and $n,m\to\infty$.
Moreover, for any $\phi\in C^\infty_0(\Omega)$, it holds that
$$
	\begin{aligned}
	& 
	\left|\int_\Omega (u_{n,m}(x,t)-u_{0}^{n,m}(x))\phi(x)\,\dee x 
	- \int_\Omega (u(x,t)-u_0(x))\phi(x)\,\dee x\right| 
	\\
	& 
	\le \|\phi\|_{L^\infty(\Omega)}(\|u_{n,m}-u\|_{L^\infty(0,\infty;L^1(\Omega))} 
	+ \|u_{0}^{n,m}-u_0\|_{L^1(\Omega)})
	\to 0\quad{\rm as}\quad n,m\to\infty,
	\end{aligned}
$$
and
$$
	\begin{aligned}
	& 
	\left|\left(\ell*\left[\int_{\Omega}\Phi(u_{n,m})\Delta\phi(x)\,\dee x\right]\right)(t) 
	- \left(\ell*\left[\int_{\Omega}\Phi(u)\Delta\phi(x)\,\dee x\right]\right)(t)\right| 
	\\
	& 
	\le \|\Delta\phi\|_{L^\infty(\Omega)}\left(\ell*\|\Phi(u_{n,m})-\Phi(u)\|_{L^1(\Omega)}\right)(t) 
	\\
	& 
	\le C\|u_{0}^{n,m}-u_0\|_{L^1(\Omega)} \to 0\ \ {\rm as}\ \ n,m\to\infty.
	\end{aligned}
$$
These imply that the function $u$ satisfies \eqref{eq:3.22}, 
namely, $u$ is an $L^1$-solution to problem \eqref{eq:IBP}.

Next, we show the uniqueness of $L^1$-solutions.
Let $u$ and $v$ be $L^1$-solutions to problem \eqref{eq:IBP} 
with initial data $u_0, v_0 \in L^1(\Omega)$, respectively.
Further, let $u_{n,m}$ and $v_{n,m}$ be approximations of $u$ and $v$, respectively.
Then, by \eqref{eq:3.27} we deduce that
$$
	\begin{aligned}
	& 
	\|u(t)-v(t)\|_{L^1_\zeta(\Omega)} 
	+\left(\ell*\|\Phi(u)-\Phi(v)\|_{L^1(\Omega)}\right)(t) 
	\\
	& 
	\le \|u_{n,m}(t)-v_{n,m}(t)\|_{L^1_\zeta(\Omega)} 
	+ \left(\ell*\|\Phi(u_{n,m})-\Phi(v_{n,m})\|_{L^1(\Omega)}\right)(t) + I_{n,m}(t) + J_{n,m}(t) 
	\\
	& 
	\le \|u_{0}^{n,m}-v_{0}^{n,m}\|_{L^1_\zeta(\Omega)} + I_{n,m}(t) + J_{n,m}(t),
	\quad t>0,
	\end{aligned}
$$
where
$$
	I_{n,m}(t) := \|u(t)-u_{n,m}(t)\|_{L^1_\zeta(\Omega)} + \|v(t)-v_{n,m}(t)\|_{L^1_\zeta(\Omega)},
$$
and
$$
	J_{n,m}(t) := \left(\ell*\|\Phi(u)-\Phi(u_{n,m})\|_{L^1(\Omega)}\right)(t) 
	+ \left(\ell*\|\Phi(v)-\Phi(v_{n,m})\|_{L^1(\Omega)}\right)(t),
$$
which converge to $0$ as $n,m\to\infty$.
This implies that \eqref{eq:3.5} holds.
Furthermore, applying an argument similar to that used in the proof of \eqref{eq:3.5},
we also obtain \eqref{eq:3.6} for $L^1$-solutions to problem \eqref{eq:IBP}.
Therefore, since $\zeta>0$ in $\Omega$, if $u_0=v_0$, then it follows that $u=v$ a.e. $x\in\Omega$ and $t>0$.
Thus the proof of Theorem~\ref{Theorem:3.2} is complete.
$\Box$
\section{Proof of main results}\label{section:4}
In this section we prove main results of this paper.
We first prove Theorem~\ref{Theorem:1.1}.

\medskip
\noindent
{\bf Proof of Theorem~\ref{Theorem:1.1}}
Let $n,m\in\N$.
We define the smooth function $\psi_n \in C_0^\infty(\R^N)$ such that
$0\le \psi_n(x)\le 1$ for $x\in\R^N$ and 
$$
	\psi_{n}(x) =
	\begin{cases}
	1 & {\rm for}\ \ 0 \le |x| \le n-1,
	\vspace{5pt}\\
	0 & {\rm for}\ \ |x| \ge n.
	\end{cases}
$$
We set
\begin{equation}
\label{eq:4.1}
	u_0^{n,m} =
	\left\{
	\begin{array}{lll}
		n\psi_n
		& {\rm if} \quad u_0>n,
		\vspace{5pt}\\
		u_0\psi_n
		& {\rm if} \quad u_0 \in [-m,n],
		\vspace{5pt}\\
		-m\psi_n
		& {\rm if} \quad u_0<-m,
	\end{array}
	\right.
\end{equation}
and consider the approximate problem
\begin{equation}
\label{Pn}\tag{$P_n$}
	\begin{cases}
		\partial_t\left(k*[u_{n,m}-u_{0}^{n,m}]\right)=\Delta\Phi(u_{n,m}), & x \in B_n,\ \ t>0,\vspace{5pt}\\
		u_{n,m}=0, & |x| \geq n,\ \ t\ge0,\vspace{5pt}\\
		u_{n,m}(x,0)=u_{0}^{n,m}(x), & x\in \R^N.
	\end{cases}
\end{equation}
It follows from Theorem~\ref{Theorem:3.2} that 
problem~\eqref{Pn} possesses a unique $L^1$-solution
$$
	u_{n,m} \in L^\infty(0,\infty;L^1(B_n))
$$
satisfying
\begin{equation*}
	\begin{aligned}
	\|u_{n,m}\|_{L^\infty(0,\infty;L^1(\R^N))} 
	&
	= \|u_{n,m}\|_{L^\infty(0,\infty; L^1(B_n))} 
	\le \|u_{0}^{n,m}\|_{L^1(B_n)}
	\\
	&
	= \|u_{0}^{n,m}\|_{L^1(\R^N)} 
	\le \|u_0\|_{L^1(\R^N)}.
	\end{aligned}
\end{equation*}
For each fixed $n\in\N$, since $u_0^{n,m} \ge u_0^{n,m+1}$ a.e. in $B_n$, it follows by Remark \ref{Remark:3.2} that
\begin{equation}
\label{eq:0}
	u_{n,m}(x,t) \ge u_{n,m+1}(x,t),\quad{\rm a.e.}\quad x\in B_n,\,\,\,t>0.
\end{equation}
Thus $\{u_{n,m}\}_m$ is monotone decreasing and uniformly bounded in $L^\infty(0,\infty;L^1(B_n))$.
Together with Theorem \ref{Theorem:3.2}, there exists an $L^1$-solution $u_n \in L^\infty(0,\infty;L^1(B_n))$ to \eqref{Pn}
such that
$$
	u_{n,m} \to u_n \quad \mbox{strongly in $L^\infty(0,\infty;L^1(B_n))$ as $m\to\infty$},
$$
and $u_n$ satisfies
\begin{equation}
\label{eq:1}
	\|u_n\|_{L^\infty(0,\infty;L^1(B_n))} = \|u_n\|_{L^\infty(0,\infty;L^1(\R^N))} \le \|u_0^n\|_{L^1(B_n)} \le \|u_0\|_{L^1(\R^N)},
\end{equation}
where $u_0^n=\min\{n,u_0\}\psi_n$.
Put
$$
	\zeta_n(x) := \left[1-\frac{|x|^2}{n^2}\right]_+.
$$
Then $\zeta_n$ satisfies
\begin{equation}
\label{eq:2}
	-\Delta\left(\frac{n^2}{2N}\zeta_n(x)\right) =
	\begin{cases}
		1 & {\rm if} \quad x\in B_n,\vspace{5pt}\\
		0 & {\rm otherwise},
	\end{cases}
\end{equation}
and $\zeta_n \to 1$, and $\Delta\zeta_n \to 0$ as $n\to\infty$.
Let $n^*\in\N$ with $n>n^*$, 
and let $u_{n^*}$ be an $L^1$-solution to approximate problem \eqref{Pn} with $n=n^*$.
Since $B_{n^*} \subset B_n$, $u_{n^*}$ can also be regarded as an $L^1$-solution to problem in $B_n$ 
with initial data $u_{0}^{n^*} \in L^1(B_{n^*}) \subset L^1(B_n)$.
Thus it follows by \eqref{eq:3.6} that
$$
	\int_{B_n}[u_{n^*}(x,t)-u_n(x,t)]_+\zeta_n(x)\,\dee x 
	\le \int_{B_n}[u_{0}^{n^*}(x)-u_{0}^{n}(x)]_+\zeta_n(x)\,\dee x,\quad t\ge0.
$$
Since $u_{0}^{n^*} \le u_{0}^{n}$ a.e. in $B_n$, we obtain
$$
	u_{n^*}(x,t) \le u_n(x,t),\quad{\rm a.e.}\quad x\in B_n,\,\,\,t>0,
$$
that is,
\begin{equation}
\label{eq:4.4}
	u_{n^*}(x,t) \le u_n(x,t),\quad{\rm a.e.}\quad x\in {\R}^N,\,\,\,t>0.
\end{equation}
This together with \eqref{eq:1} implies that
$\{u_n\}_n$ is monotone increasing and uniformly bounded in $L^\infty(0,\infty;L^1(\R^N))$.
Hence there exists a $u \in L^\infty(0,\infty;L^1(\R^N))$ such that
$$
	u_n \to u\quad\mbox{strongly in $L^\infty(0,\infty;L^1({\R}^N))$ as $n\to\infty$}.
$$
Furthermore, we see that
\begin{equation}
\label{eq:01}
	u_n\zeta_n \to u\quad\mbox{strongly in $L^\infty(0,\infty;L^1({\R}^N))$ as $n\to\infty$}.
\end{equation}
In fact, we have
$$
	\begin{aligned}
	\|u_n\zeta_n-u\|_{L^\infty(0,\infty;L^1(\R^N))}
	& 
	\le \|u-\zeta_nu\|_{L^\infty(0,\infty;L^1(\R^N))} 
	+ \|u_n-u\|_{L^\infty(0,\infty;L^1_{\zeta_n}(\R^N))} \\
	& =: I_n +J_n.
	\end{aligned}
$$
For the first term in the right-hand side, it holds that
$$
	I_n
	\le \|u\|_{L^\infty(0,\infty;L^1(\R^N))}\|1-\zeta_n\|_{L^\infty(\R^N)}
	= \|u\|_{L^\infty(0,\infty;L^1(\R^N))}
	< \infty.
$$
Since $\zeta_n\to1$ as $n\to\infty$, 
applying the dominated convergence theorem, we obtain
$I_n\to0$ as $n\to\infty$.
On the other hand, since $|\zeta_n| \le 1$, we have
$$
	J_n
	\le \|u_n-u\|_{L^\infty(0,\infty;L^1(\R^N))}
	\to 0\quad{\rm as}\quad n\to\infty.
$$
These yield \eqref{eq:01}.
On the other hand, by Theorem~\ref{Theorem:3.2}, 
we observe that $u_n\in C([0,\infty);L^1_{\zeta_n}(B_n))$, that is, $u_n \in C([0,\infty);L^1_{\zeta_n}({\R}^N))$.
Then, for $t,s\ge0$, we deduce that
$$
\begin{aligned}
	&
	\|u(t)-u(s)\|_{L^1(\R^N)}
	\\
	& 
	\le \|u(t)-u_n(t)\zeta_n\|_{L^1(\R^N)} 
	+ \|u(s)-u_n(s)\zeta_n\|_{L^1(\R^N)} 
	+ \|u_n(t)-u_n(s)\|_{L^1_{\zeta_n}(\R^N)} 
	\\
	&
	\le 2\|u-u_n\zeta_n\|_{L^\infty(0,\infty;L^1(\R^N))} 
	+ \|u_n(t)-u_n(s)\|_{L^1_{\zeta_n}(\R^N)}.
\end{aligned}
$$
This together with \eqref{eq:1} and \eqref{eq:01} implies that  $u \in C([0,\infty);L^1(\R^N))$ 
and $\|u(t)\|_{L^1(\R^N)}$ is bounded on $[0,\infty)$.

From now, we show local integrability on $\R^N$ of $\Phi(u)$.
Let $M>0$ be an arbitrary constant.
We now choose sufficiently large $n\in\N$ satisfying $B_M \subset B_n$.
Let $m^*\in\N$, and put
$$
	U:=u_{n,m}-u_{n,m^*}, \quad U_0:=u_{0}^{n,m}-u_{0}^{n,m^*}, \quad V:=\Phi(u_{n,m})-\Phi(u_{n,m^*}).
$$
Here we recall that $u_{0}^{n,m} \in L^\infty(B_n)$, 
and $u_{n,m}$ is also a global weak solution to problem \eqref{Pn}.
Hence $u_{n,m}$ and $u_{n,m^*}$
satisfy the time-regularized version of the weak formulation of problem \eqref{Pn}, 
that is, for any $\varphi_n \in H_0^1(B_n)$,
\begin{equation}
\label{eq:5}
	\begin{aligned}
		&
		\int_{B_n}\partial_t(k_j*[U(x,\cdot)-U_0(x)])(t)\varphi_n(x)\, \dee x 
		\\
		& \hspace{2cm}
		+ \int_{B_n}\left(h_j*\left[\nabla V(x,\cdot)\right]\right)(t)\cdot\nabla\varphi_n(x)\,\dee x
		=0, 
		\quad{\rm a.e.}\,\,\, t>0, \quad j\in\N,
	\end{aligned}
\end{equation}
where $h_j$ is given in \eqref{eq:2.3} and $k_j = k*h_j$.
Let $H_{\E}$ be function given in \eqref{eq:3.19} with $p=1$.
Then, taking the function $\varphi=H_{\E}'(V)\zeta_M \in H_0^1(B_n)$ as a test function in \eqref{eq:5},
and applying the same argument as in the proof of \eqref{eq:11}, we see that
\begin{equation}
\label{eq:7}
	\begin{aligned}
	&
	\int_{B_n}\bigg(H_{\E}(U(x,t))-H_{\E}(U_0(x))\bigg)\zeta_M(x)\,\dee x
	\\
	&\qquad
	+ \left(\ell*\left[\int_{B_n}|\nabla V(x,\cdot)|^2H_{\E}''(V(x,\cdot))\zeta_M(x)\,\dee x\right]\right)(t) 
	\\
	&\qquad\qquad
	+
	\left(\ell*\left[\int_{B_n}\nabla H_{\E}(V(x,\cdot))\cdot\nabla\zeta_M(x)\,\dee x\right]\right)(t)
	\le |(\ell*J_{\E})|(t),
	\quad{\rm a.e.}\,\,\, t>0,
	\end{aligned}
\end{equation}
where
$$
	J_{\E}(t) 
	= \int_{B_n}\partial_t(k*[U(x,\cdot)-U_0(x)])(t)\biggl\{H_{\E}'(U(x,t))-H_{\E}'(V(x,t))\biggl\}\zeta_M(x)\,\dee x.
$$
Since $\zeta_M\ge0$ in $B_n$, the second term in the left-hand side of \eqref{eq:7} is nonnegative.
Furthermore, for the third term in the left-hand side, by \eqref{eq:2}, we obtain
$$
	\begin{aligned}
	&
	\int_{B_n}\nabla H_{\E}(V(x,t))\cdot\nabla\zeta_M(x)\,\dee x
	= -\int_{B_n}H_{\E}(V(x,t))\Delta\zeta_M(x)\,\dee x \\
	&\qquad\qquad
	= \frac{2N}{M^2}\int_{B_n}H_{\E}(V(x,t))\chi_{B_M}\,\dee x
	= \frac{2N}{M^2}\int_{B_M}H_{\E}(V(x,t))\,\dee x,
	\quad{\rm a.e.}\,\,\,t>0,
	\end{aligned}
$$
and hence deduce that
$$
	\begin{aligned}
		&
		\int_{B_M}H_{\E}(U(x,t))\zeta_M(x)\,\dee x
		+ \frac{2N}{M^2}\left(\ell*\left[\int_{B_M}H_{\E}(V(x,\cdot))\,\dee x\right]\right)(t) \\
		& \qquad\qquad
		\le \int_{B_M}H_{\E}(U_0(x))\zeta_M(x)\,\dee x + |(\ell*J_{\E})|(t),
		\quad{\rm a.e.}\,\,\, t>0.
	\end{aligned}
$$
Moreover, applying the same argument as in the proof of \eqref{eq:12}, 
it holds that $|(\ell*J_{\E})|(t)\to 0$ as $\E \to 0^+$.
Since $H_{\E}(y)\to|y|$ as $\E\to0^+$ for all $y\in\R$,
applying the dominated convergence theorem, we have
$$
	\begin{aligned}
		&
		\|u_{n,m}(t)-u_{n,m^*}(t)\|_{L^1_{\zeta_M}(B_M)}
		+ \frac{2N}{M^2}\left(\ell*\|\Phi(u_{n,m})-\Phi(u_{n,m^*})\|_{L^1(B_M)}\right)(t) \\
		& \hspace{5cm}
		\le \|u_{0}^{n,m}-u_{0}^{n,m^*}\|_{L^1_{\zeta_M}(B_M)},
		\quad{\rm a.e.}\,\,\, t>0.
	\end{aligned}
$$
Furthermore, by $|\zeta_M| \le 1$ and $B_M\subset B_n$, we see that
$$
	\|u_{0}^{n,m}-u_{0}^{n,m^*}\|_{L^1_{\zeta_M}(B_M)}
	\le \|u_{0}^{n,m}-u_{0}^{n,m^*}\|_{L^1(B_n)}.
$$
Since $u_{n,m} \in C([0,\infty);H^{-1}(B_n))$, $\zeta_M \in H_0^1(B_n)$,
it holds that $u_{n,m}\in C([0,\infty);L^1_{\zeta_M}(B_n))$.
Then, by Theorem~\ref{Theorem:3.2}, we have $\ell*\|\Phi(u_{n,m})\|_{L^1(B_n)} \in C([0,\infty))$.
Therefore, it follows that, for any $M>0$,
\begin{equation}
\label{eq:4-1}
	\begin{aligned}
		&
		\|u_{n,m}(t)-u_{n,m^*}(t)\|_{L^1_{\zeta_M}(B_M)} 
		+ \frac{2N}{M^2}\left(\ell*\|\Phi(u_{n,m})-\Phi(u_{n,m^*})\|_{L^1(B_M)}\right)(t) \\
		&\hspace{5cm}
		\le \|u_{0}^{n,m}-u_{0}^{n,m^*}\|_{L^1(B_n)}, \quad t\ge0.
	\end{aligned}
\end{equation}
For any fixed $T>0$, convolving this with the kernel $k$, and using $k*\ell=1$, we also get
\begin{equation}
\label{eq:4-2}
	\frac{2N}{M^2}\|\Phi(u_{n,m})-\Phi(u_{n,m^*})\|_{L^1(0,T;L^1(B_n))} \le \|k\|_{L^1(0,T)}\|u_{0}^{n,m}-u_{0}^{n,m^*}\|_{L^1(B_n)}.
\end{equation}
This together with Theorem \ref{Theorem:3.2} and condition ($L$), there exists a $u_n$ such that
$$
	\begin{aligned}
	u_{n,m} \to u_n & \quad \mbox{strongly in $C([0,\infty); L_{\zeta_M}^1(B_n))$ as $m\to\infty$}, \\
	\Phi(u_{n,m}) \to \Phi(u_n) & \quad \mbox{strongly in $L_{\rm loc}^1(0,\infty;L^1(B_n))$ as $m\to\infty$}, \\
	\ell*\|\Phi(u_{n,m})\|_{L^1(B_n)} \to \ell*\|\Phi(u_n)\|_{L^1(B_n)} & \quad \mbox{strongly in $C([0,\infty))$ as $m\to\infty$},
	\end{aligned}
$$
and it holds that
\begin{equation}
\label{eq:4-2.5}
	\left(\ell*\|\Phi(u_n)\|_{L^1(B_n)}\right)(0)=0.
\end{equation}
Moreover, let $n^*\in\N$ satisfying $B_M \subset B_{n^*} \subset B_n$, by \eqref{eq:4-1} and \eqref{eq:4-2}, $u_n$ and $u_{n^*}$  satisfy
\begin{equation}
\label{eq:4-3}
	\begin{aligned}
	\|u_n(t)-u_{n^*}(t)\|_{L_{\zeta_M}^1(B_M)} + \frac{2N}{M^2}\left(\ell*\|\Phi(u_n)-\Phi(u_{n^*})\|_{L^1(B_M)}\right)(t) \\
	\le\|u_0^n-u_0^{n^*}\|_{L^1(B_n)},\,\,\, t\ge0,
	\end{aligned}
\end{equation}
and
\begin{equation*}
	\frac{2N}{M^2}\|\Phi(u_n)-\Phi(u_{n^*})\|_{L^1(0,T;L^1(B_M))}
	\le \|k\|_{L^1(0,T)}\|u_0^n-u_0^{n^*}\|_{L^1(B_n)}
\end{equation*}
for any fixed $T>0$.
By \eqref{eq:4-3}, we deduce that
$$
	\begin{aligned}
		&
		\sup_{t\ge0}\bigg|\left(\ell*\|\Phi(u_{n})\|_{L^1(B_M)}\right)(t) 
		- \left(\ell*\|\Phi(u_{n^*})\|_{L^1(B_M)}\right)(t)\bigg| \\
		&
		\le \sup_{t\ge0}\bigg|\left(\ell*\|\Phi(u_{n})-\Phi(u_{n^*})\|_{L^1(B_M)}\right)(t)\bigg|
		\le \|u_{0}^{n}-u_{0}^{n^*}\|_{L^1(B_n)}.
	\end{aligned}
$$
Thus, by arbitrariness of $M$ and $T$, 
we see that there exist $w \in L_{\rm loc}^1(0,\infty;L_{\rm loc}^1(\R^N))$ and $z \in C([0,\infty))$ such that
$$
	\begin{aligned}
		\Phi(u_n) & \to w \quad \mbox{strongly in $L_{\rm loc}^1(0,\infty;L_{\rm loc}^1(\R^N))$ as $n\to\infty$}, 
		\\
		\ell*\|\Phi(u_n)\|_{L_{\rm loc}^1(\R^N)} & \to z \quad \mbox{strongly in $C([0,\infty))$ as $n\to\infty$}.
	\end{aligned}
$$
Since $u_n \to u$ a.e. by choosing appropriate subsequence, 
we observe that $w=\Phi(u)$ a.e., and $z=\ell*\|\Phi(u)\|_{L_{\rm loc}^1(\R^N)}$.
Moreover, it also holds that
$$
	\left(\ell*\|\Phi(u)\|_{L_{\rm loc}^1(\R^N)}\right)(0)=0.
$$
Indeed, by \eqref{eq:4-2.5} and \eqref{eq:4-3}, we have
\begin{align*}
	\left(\ell*\|\Phi(u)\|_{L_{\rm loc}^1(\R^N)}\right)(t)
	&
	\le \left(\ell*\|\Phi(u)-\Phi(u_n)\|_{L_{\rm loc}^1(\R^N)}\right)(t)
	+ \left(\ell*\|\Phi(u_n)\|_{L_{\rm loc}^1(\R^N)}\right)(t) \\
	&
	\le C\|u_0-u_{0}^n\|_{L^1(\R^N)}
	+ \left(\ell*\|\Phi(u_n)\|_{L_{\rm loc}^1(\R^N)}\right)(t)
	\to 0
\end{align*}
as $t\to0^+$ and $n\to\infty$.

Finally we prove that $L^1$-contraction principle holds.
Let $u_n$ and $v_n$ be approximations of $u$ and $v$, respectively. 
Then, by \eqref{eq:3.5} with $\zeta_n$ to $u_n$ and $v_n$, we have
$$
	\|u_n(t)-v_n(t)\|_{L^1_{\zeta_n}(B_n)}
	\le \|u_{0,n}-v_{0,n}\|_{L^1_{\zeta_n}(B_n)}, \quad t\ge0.
$$
This implies that
$$
	\begin{aligned}
	&
	\|u(t)-v(t)\|_{L^1(\R^N)}
	\\
	&
	\le \|u(t)-u_n(t)\zeta_n\|_{L^1(\R^N)}
	+ \|v(t)-v_n(t)\zeta_n\|_{L^1(\R^N)}
	+ \|u_n(t)-v_n(t)\|_{L^1_{\zeta_n}(B_n)} 
	\\
	&
	\le \|u-u_n\zeta_n\|_{L^\infty(0,\infty;L^1(\R^N))}
	+ \|v-v_n\zeta_n\|_{L^\infty(0,\infty;L^1(\R^N))}
	+ \|u_{0,n}-v_{0,n}\|_{L^1_{\zeta_n}(B_n)} 
	\\
	&
	\to \|u_0-v_0\|_{L^1(\R^N)} \quad \mbox{as $n\to\infty$}.
	\end{aligned}
$$
Thus \eqref{eq:1.4} holds.
Furthermore, applying an argument similar to that used in the proof of \eqref{eq:1.4}, we obtain \eqref{eq:1.3}.
Therefore, uniqueness of $L^1$-solutions to problem~\eqref{eq:P} holds,
and the proof of Theorem~\ref{Theorem:1.1} is complete. 
$\Box$
\bigskip

Similarly to Theorem~\ref{Theorem:3.1},
we show the $L^p$-estimate for $L^1$-solutions to problem \eqref{eq:P}.
\begin{proposition}
\label{Proposition:4.1}
	Assume the same conditions as in Theorem~$\ref{Theorem:1.1}$.
	Furthermore, let $u_0 \in L^p(\R^N)$ for some $p>1$.
	Let $u \in C([0,\infty);L^1(\R^N))$ be an $L^1$-solution to problem \eqref{eq:P}.
	Then $u$ belongs to $L^\infty(0,\infty;L^p(\R^N))$ and satisfies
	$$
	\|u\|_{L^\infty(0,\infty;L^p(\R^N))} \leq \|u_0\|_{L^p(\R^N)}.
	$$
\end{proposition}
{\bf Proof.}
Let $u_{n,m}$ be an $L^1$-solution to approximate problem \eqref{Pn} 
with initial data $u_{0}^{n,m}$ which is defined as \eqref{eq:4.1}.
Then, by Theorem~\ref{Theorem:3.1}, we have
$$
	\|u_{n,m}\|_{L^\infty(0,\infty;L^p(\R^N))}
	= \|u_{n,m}\|_{L^\infty(0,\infty;L^p(B_n(0)))}
	\le \|u_{0}^{n,m}\|_{L^p(\R^N)}
	\le\|u_0^{n,m}\|_{L^p(\R^N)}.
$$
Furthermore, by \eqref{eq:0} we see that 
$\{u_{n,m}\}_m$ is monotone decreasing and uniformly bounded in $L^\infty(0,\infty;L^p(B_n))$.
Therefore, there exists a $u_n$ such that
$$
	u_{n,m} \to u_n \quad \mbox{strongly in $L^\infty(0,\infty;L^p(B_n))$ as $m\to\infty$},
$$
and $u_n$ satisfies
$$
	\|u_n\|_{L^\infty(0,\infty;L^p(\R^N))}
	= \|u_n\|_{L^\infty(0,\infty;L^p(B_n))}
	\le \|u_0^n\|_{L^p(B_n)}
	\le \|u_0\|_{L^p(\R^N)}.
$$
Furthermore, from \eqref{eq:4.4}, we see that $\{u_n\}_n$ is monotone increasing and uniformly bounded in $L^\infty(0,\infty;L^p(\R^N))$.
Hence there exists a $u$ such that
$$
	u_n \to u \quad \mbox{strongly in $L^\infty(0,\infty;L^p(\R^N))$ as $n\to\infty$}.
$$
Combining with Theorem \ref{Theorem:1.1}, $u_n$ also converges in $C([0,\infty);L^1(\R^N))$.
Thus Proposition~\ref{Proposition:4.1} follows.
$\Box$
\bigskip

Next we prove Theorem~\ref{Theorem:1.2}.

\medskip
\noindent
{\bf Proof of Theorem~\ref{Theorem:1.2}.}
We use the test function introduced in \cite{HP}.
We define $\varphi \in C_0^\infty(\R^N)$ such that
$$
	\varphi(x)=
	\begin{cases}
	1 & {\rm if}\ \ |x| \le 1,\vspace{5pt}\\
	0 & {\rm if}\ \ |x| \ge 2,
\end{cases}
\hspace{1cm} 0\le \varphi \le 1,
$$
and set $\varphi_R(x)=\varphi\left(x/R\right)$.
Note that $\Delta\varphi_R(x) = R^{-2}\Delta\varphi\left(x/R\right)$.
Then by choosing $\phi=\varphi_R$ as a test function 
and using H\"{o}lder's inequality in the right-hand side in \eqref{eq:1.1}, we obtain
$$
	\left|\int_{B_{2R}}\left(u(x,t)-u_0(x)\right)\,\dee x\right| 
	\le I(\varphi_R)\left(\ell*\left[\left(\int_{B_{2R}\setminus B_R}|u(x,\cdot)|\,\dee x\right)^m\right]\right)(t),
$$
where
$$
	I(\varphi_R)=\left(\int_{B_{2R}\setminus B_R}\left|\Delta\varphi_R(x)\right|^{\frac{1}{1-m}}\,\dee x\right)^{1-m}.
$$
Then we have
$$
	I(\varphi_R) \le C_1R^{N(1-m)-2} = C_1R^{-N\left(m-\frac{N-2}{N}\right)},
$$
where $C_1$ is a positive constant which depends only on $m$ and $N$.
Since $u \in L^\infty(0,\infty;L^1(\R^N))$, it follows that
$$
	\left|\int_{B_{2R}}(u(x,t)-u_0(x))\ dx\right| 
	\le C_1\|\ell\|_{L^1(0,t)}
	\sup_{t>0}\left(\int_{B_{2R}\setminus B_R}|u(x,t)|\,\dee x\right)^mR^{-N\left(m-\frac{N-2}{N}\right)} \to 0
$$
as $R\to\infty$ for all $t>0$.
This implies the desired identity \eqref{eq:1.5},
thus the proof of Theorem~\ref{Theorem:1.2} is complete.
$\Box$
\medskip

\begin{remark}
\label{Remark:4.1}
	If $\Phi(r)=|r|^{m-1}r$ with $m>1$, 
	the identity \eqref{eq:1.5} is shown by assuming $u_0 \in L^1(\R^N) \cap L^p(\R^N)$ with $p\ge m$.
	Indeed, let $\psi_{n}$ be given in the proof of Theorem~$\ref{Theorem:1.1}$. 
	Taking $\phi=\psi_{n}$ as a test function in the formulation \eqref{eq:1.1}, we obtain
	\begin{equation}
	\label{eq:4.6}
		\left|\int_{B_n}\left(u(x,t)-u_0(x)\right)\psi_n(x)\,\dee x\right|
		\le \left(\ell*\left[\int_{B_n\setminus B_{n-1}}|u(x,\cdot)|^m|\Delta\psi_{n}(x)|\,\dee x\right]\right)(t).
	\end{equation}
	For the left-hand side, it holds that
	$$
		\left|\int_{B_n}\left(u(x,t)-u_0(x)\right)\psi_{n}(x)\,\dee x\right|
		\to \left|\int_{\R^N}\left(u(x,t)-u_0(x)\right)\,\dee x\right|\quad\mbox{as $n\to\infty$},
	$$
	for all $t>0$.
	On the other hand, for the right-hand side of \eqref{eq:4.6}, 
	since $u \in L^\infty(0,\infty; L^q(\R^N))$ for $q\in[1,p]$,
	we see that the right-hand side tends to $0$ as $n\to\infty$.
	Consequently, we conclude that
	$$
		\int_{\R^N}u(x,t)\,\dee x = \int_{\R^N}u_0(x)\,\dee x,\quad t>0.
	$$
\end{remark}
\bigskip

Finally we prove Theorem~\ref{Theorem:1.3}.
Before proving Theorem~\ref{Theorem:1.3},
we recall a result on asymptotic behavior of a nonlinear fractional ODE.
\begin{lemma}\cite{VZ2}*{Theorem 7.1}
\label{Lemma:4.1}
	Let $\alpha \in (0,1)$, $\lambda > 0$, and $m>0$.
	Assume that $v \in W_{\rm loc}^{1,1}([0,\infty))$ is the solution to
	$$
		\RL\left(v-v_0\right) = -\lambda v^m,\qquad v(0)=v_0.
	$$
	Then there exist $c_1,c_2>0$ depending only on $v_0$, $m$, $\alpha$, and $\lambda$ such that
	$$
		\frac{c_1}{1+t^{\frac{\alpha}{m}}} \leq v(t) \leq \frac{c_2}{1+t^{\frac{\alpha}{m}}},\ \ t>0.
	$$
\end{lemma}
\bigskip

Now we are ready to prove Theorem~\ref{Theorem:1.3}.

\noindent
{\bf Proof of Theorem~\ref{Theorem:1.3}.}
This proof is very similar to that used in \cite{BGI}*{Proposition 5.2} after choosing an appropriate test function.
However, for the reader's convenience, we present the details.

Let $\psi_n$ be given in the proof of Theorem \ref{Theorem:1.1}. 
Put $G(x)=\exp({-|x|^2})$ for $x\in\R^N$.
Then we see that $G \in C^\infty(\R^N)$ is the positive smooth function satisfying
$G, \nabla G, \Delta G \in L^q({\R}^N)$ 
for all $q\in[1,\infty]$
and
\begin{equation}
\label{eq:4.7}
	\Delta G(x) \ge -2N G(x)\quad{\rm in}\,\,\,{\R}^N.
\end{equation}
Then, taking $\phi(x)=\psi_n(x)G(x)$ as a test function in the formulation \eqref{eq:1.1}, we obtain
$$
	\int_{\R^N}\left(u(x,t)-u_0(x)\right)\psi_n(x)G(x)\,\dee x 
	= \biggl(\ell*\left[J_1 + J_2 + J_3\right]\biggr)(t),\,\,\, t>0,
$$
where
$$
	\begin{aligned}
	J_1(t) &:= \int_{\R^N}u^m(x,t)G(x)\Delta\psi_n(x)\,\dee x,
	\\
	J_2(t) &: = 2\int_{\R^N}u^m(x,t)\nabla G(x)\cdot\nabla\psi_n(x)\,\dee x,
	\\
	J_3(t) &:= \int_{\R^N}u^m(x,t)\psi_n(x)\Delta G(x)\,\dee x.
	\end{aligned}
$$
Here we note that $J_1,J_2 \to 0$ as $n\to\infty$.
Indeed, since $u \in L^\infty(0,\infty;L^1(\R^N))$, we have
$$
	\begin{aligned}
	|J_1(t)|
	& 
	= \left|\int_{B_n\setminus B_{n-1}(0)}u^m(x,t)G(x)\Delta\psi_n(x)\,\dee x\right| 
	\\
	& 
	\le \|\Delta\psi_n\|_{L^\infty(\R^N)}\left(\int_{B_n\setminus B_{n-1}}u(x,t)\, \dee x\right)^m
	\left(\int_{B_n\setminus B_{n-1}}G^{\frac{1}{1-m}}(x)\,\dee x\right)^{1-m} ,
	\end{aligned}
$$
and
$$
	\begin{aligned}
	|J_2(t)|
	& 
	= \left|2\int_{B_n\setminus B_{n-1}}u^m(x,t)\nabla G(x)\cdot\nabla\psi_n(x)\,\dee x\right| 
	\\
	& 
	\le 2\|\nabla\psi_n\|_{L^\infty(\R^N)}\left(\int_{B_n\setminus B_{n-1}}u(x,t)\, \dee x\right)^m
	\left(\int_{B_n\setminus B_{n-1}}|\nabla G(x)|^{\frac{1}{1-m}}\,\dee x\right)^{1-m}
	\end{aligned}
$$
for $t>0$. 
Thus, letting $n\to\infty$, we obtain
\begin{equation}
\label{eq:4.9}
	\int_{\R^N}\left(u(x,t)-u_0(x)\right)G(x)\,\dee x 
	= \left(\ell*\left[\int_{\R^N}u^m(x,\cdot)\Delta G(x)\,\dee x\right]\right)(t),\quad t>0.
\end{equation}
Here since $m\in(0,1)$, we see that each term is well-defined for each $t>0$.
Convolving \eqref{eq:4.9} with $k_j \in W_{\rm loc}^{1,1}([0,\infty))$ and using $k*\ell = 1$, we obtain
\begin{equation}
\label{eq:4.10}
	\begin{aligned}
	&
	\left(k_j*\left[\int_{\R^N}\left(u(x,\cdot)-u_0(x)\right)G(x)\,\dee x\right]\right)(t)
	\\
	& 
	= \int_0^t\left(h_j*\left[\int_{\R^N}u^m(x,\cdot)\Delta G(x)\,\dee x\right]\right)(s)\,\dee s,
	\quad t>0.
	\end{aligned}
\end{equation}
Since $k_j \in W_{\rm loc}^{1,1}([0,\infty))$ 
and $u \in C([0,\infty);L^1(\R^N)) \cap L^\infty(0,\infty;L^1(\R^N))$, we see that
$$
	\left(k_j*\left[\int_{\R^N}\left(u(x,\cdot)-u_0(x)\right)G(x)\,\dee x\right]\right) \in W_{\rm loc}^{1,1}([0,\infty)).
$$
Thus both sides of \eqref{eq:4.10} are differentiable for a.e. $t>0$, so we get
$$	\begin{aligned}
	&
	\dt\left(k_j*\left[\int_{\R^N}\left(u(x,\cdot)-u_0(x)\right)G(x)\,\dee x\right]\right)(t) 
	\\
	&
	= \left(h_j*\left[\int_{\R^N}u^m(x,\cdot)\Delta G(x)\,\dee x\right]\right)(t),
	\quad {\rm a.e.}\,\,\, t>0.
	\end{aligned}
$$
Since it follows from H\"{o}lder's inequality with \eqref{eq:4.7} that
$$
	\begin{aligned}
	\int_{\R^N}u^m(x,t)\Delta G(x)\,\dee x
	& 
	\ge -2N\int_{\R^N}u^m(x,t)G(x)\,\dee x 
	\\
	& 
	\ge -2N\left(\int_{\R^N}G(x)\ dx\right)^{1-m}\left(\int_{\R^N}u(x,t)G(x)\,\dee x\right)^m,
	\quad t>0,
	\end{aligned}
$$
we have
\begin{equation}
\label{eq:4.11}
	\dt\left(k_j*[U(\cdot)-U_0]\right)(t) \ge -C\left(h_j*U^m(\cdot)\right)(t),\quad{\rm a.e.}\,\,\,t>0,
\end{equation}
where $C$ depends only on $N$ and $m$,
and
$$
	U(t) := \int_{\R^N}u(x,t)G(x)\,\dee x,\qquad
	U_0 = \int_{\R^N}u_0(x)G(x)\,\dee x.
$$
Let $Z \in W_{\rm loc}^{1,1}([0,\infty))$ be a solution to the Cauchy problem
$$
		\RL\left(Z(\cdot)-Z_0\right)(t) = -CZ^m(t), \qquad
		Z(0)=U_0.
$$
Convolving with $h_j$ defined as \eqref{eq:2.3}, and taking the difference with \eqref{eq:4.11}, we obtain
\begin{equation}
\label{eq:4.12}
	\dt\left(k_j*[Z(\cdot)-U(\cdot)]\right)(t) + C\left(h_j*[Z^m(\cdot)-U^m(\cdot)]\right)(t) 
	\le 0,\quad{\rm a.e.}\,\,\,t>0.
\end{equation}
In order to prove that $Z(t) \le U(t)$ for a.e. $t>0$, we assume that $Z(t)>U(t)$ for a.e. $t>0$.
Multiplying \eqref{eq:4.12} by $[Z-U]_+ \in L^\infty([0,\infty))$,
and applying Lemma \ref{Lemma:2.1} for $H(y) = \frac{1}{2}[y]_+^2$, we see that
$$
	\frac{1}{2}\dt\left(k_j*[Z(\cdot)-U(\cdot)]_+^2\right)(t) 
	+ C\left(h_j*[Z^m(\cdot)-U^m(\cdot)]\right)(t)[Z(t)-U(t)]_+ \le 0,
	\quad{\rm a.e.}\,\,\, t>0.
$$
Since $\left(k_j*[Z(\cdot)-U(\cdot)]_+^2\right)(0)=0$, 
for any $T>0$,
by integrating both sides from $0$ to $T$, we have
\begin{equation}
\label{eq:4.13}
	\frac{1}{2}\left(k_j*[Z(\cdot)-U(\cdot)]_+^2\right)(T) 
	+ C\int_0^T\left(h_j*[Z^m(\cdot)-U^m(\cdot)]\right)(t)[Z(t)-U(t)]_+\,\dee t \le 0.
\end{equation}
Furthermore, it holds that
$$
	\begin{aligned}
	& 
	\left|\int_0^T\left(h_j*[Z^m(\cdot)-U^m(\cdot)]\right)(t)[Z(t)-U(t)]_+\,\dee t 
	- \int_0^T\left(Z^m(t)-U^m(t)\right)[Z(t)-U(t)]_+\,\dee t\right| 
	\\
	& 
	\le \|[Z-U]_+\|_{L^\infty(0,T)}\left\|\left(h_j*[Z^m-U^m]\right)-(Z^m-U^m)\right\|_{L^1(0,T)} 
	\\
	& 
	\to 0\quad{\rm as}\,\,\, j\to\infty,
	\end{aligned}
$$
and
$$
	\begin{aligned}
	&
	\left|\left(k_j*[Z(\cdot)-U(\cdot)]_+^2\right)(T)-\left(g_\alpha*[Z(\cdot)-U(\cdot)]_+^2\right)(T)\right|
	\\
	& \le \|[Z-U]_+\|_{L^\infty(0,T)}^2\|k_j-g_\alpha\|_{L^1(0,T)}
	\to 0\quad{\rm as}\,\,\, j\to\infty.
	\end{aligned}
$$
Thus, taking the limit $j\to\infty$ in \eqref{eq:4.13}, we obtain
$$
	\frac{1}{2}\left(g_\alpha*[Z(\cdot)-U(\cdot)]_+^2\right)(T) 
	+ \int_0^T\left(Z^m(t)-U^m(t)\right)[Z(t)-U(t)]_+\,\dee t \le 0.
$$
Since $\left(Z^m-U^m\right)[Z-U]_+ \ge 0$ and $g_\alpha$ is nonnegative,
by the arbitrariness of $T$ we see that
$$
	Z(t) \leq U(t),\quad{\rm a.e.}\,\,\,t>0.
$$
Therefore, it follows from Lemma~\ref{Lemma:4.1} that
$$
	\frac{c_1}{1+t^{\frac{\alpha}{m}}} \le Z(t) \le U(t),\quad t>0.
$$
On the other hand, applying H\"{o}lder's inequality, we have
$$
	U(t)
	=\int_{\R^N}u(x,t)G(x)\,\dee x
	\le \|u(t)\|_{L^p(\R^N)}\|G\|_{L^{q}(\R^N)} 
$$
for $p\in[1,\infty]$ with $1/p+1/q=1$.
Thus we obtain \eqref{eq:1.6}, and the proof of Theorem~\ref{Theorem:1.3} is complete.
$\Box$
\bigskip

\noindent
{\bf Acknowledgements.} 
T. K. was supported in part by JSPS KAKENHI Grant Number JP25H00591 and JP22KK0035.



\begin{bibdiv}
\begin{biblist}
\bib{A}{article}{
    author={Akagi, Goro},
    title={Fractional flows driven by subdifferentials in {H}ilbert spaces},
    journal={Israel J. Math.},
    volume={234},
    date={2019},
    number={2},
    pages={809--862},
     issn={0021-2172,1565-8511},
}
\bib{AN}{article}{
      title={Time-fractional gradient flows for nonconvex energies in Hilbert spaces}, 
      author={Akagi, Garo},
      author={Nakajima, Yoshihito},
      year={2025},
      eprint={2501.08059},
}
\bib{AAC}{article}{
   author={Achleitner, Franz},
   author={Akagi, Goro},
   author={Kuehn, Christian},
   author={Melenk, Jens Markus},
   author={Rademacher, Jens D. M.},
   author={Soresina, Cinzia},
   author={Yang, Jichen},
   title={Fractional dissipative PDEs},
   conference={
      title={Fractional dispersive models and applications---recent
      developments and future perspectives},
   },
   book={
      series={Nonlinear Syst. Complex.},
      volume={37},
      publisher={Springer, Cham},
   },
   isbn={978-3-031-54977-9},
   isbn={978-3-031-54978-6},
   date={[2024] \copyright 2024},
   pages={53--122},
}
\bib{AB}{article}{
   author={Alves, Claudianor O.},
   author={Boudjeriou, Tahir},
   title={Existence and uniqueness of solution for some time fractional
   parabolic equations involving the 1-Laplace operator},
   journal={Partial Differ. Equ. Appl.},
   volume={4},
   date={2023},
   number={2},
   pages={Paper No. 5, 35},
   issn={2662-2963},
}
\bib{B}{book}{
   author={Barbu, Viorel},
   title={Nonlinear semigroups and differential equations in Banach spaces},
   note={Translated from the Romanian},
   publisher={Editura Academiei Republicii Socialiste Rom\^ania, Bucharest;
   Noordhoff International Publishing, Leiden},
   date={1976},
   pages={352},
}
\bib{BC}{article}{
   author={Br\'ezis, Ha\"im},
   author={Crandall, Michael G.},
   title={Uniqueness of solutions of the initial-value problem for
   $u\sb{t}-\Delta \varphi (u)=0$},
   journal={J. Math. Pures Appl. (9)},
   volume={58},
   date={1979},
   number={2},
   pages={153--163},
   issn={0021-7824},
}
\bib{BGI}{article}{
   author={Bonforte, Matteo},
   author={Gualdani, Maria},
   author={Ibarrondo, Peio},
   title={Time-fractional porous medium type equations. sharp time decay and
   regularization},
   journal={Calc. Var. Partial Differential Equations},
   volume={65},
   date={2026},
   number={1},
   pages={Paper No. 30},
   issn={0944-2669},
}
\bib{C}{article}{
   author={Cl\'ement, Ph.},
   title={On abstract Volterra equations in Banach spaces with completely
   positive kernels},
   conference={
      title={Infinite-dimensional systems},
      address={Retzhof},
      date={1983},
   },
   book={
      series={Lecture Notes in Math.},
      volume={1076},
      publisher={Springer, Berlin},
   },
   isbn={3-540-13376-3},
   date={1984},
   pages={32--40},
}
\bib{CQW}{article}{
   author={Cort\'azar, Carmen},
   author={Quir\'os, Fernando},
   author={Wolanski, Noem\'i},
   title={A heat equation with memory: large-time behavior},
   journal={J. Funct. Anal.},
   volume={281},
   date={2021},
   number={9},
   pages={Paper No. 109174, 40},
   issn={0022-1236},
}
\bib{CP}{article}{
   author={Crandall, Michael G.},
   author={Pierre, Michel},
   title={Regularizing effects for $u\sb{t}=\Delta \varphi (u)$},
   journal={Trans. Amer. Math. Soc.},
   volume={274},
   date={1982},
   number={1},
   pages={159--168},
   issn={0002-9947},
}
\bib{EK}{article}{
   author={Eidelman, Samuil D.},
   author={Kochubei, Anatoly N.},
   title={Cauchy problem for fractional diffusion equations},
   journal={J. Differential Equations},
   volume={199},
   date={2004},
   number={2},
   pages={211--255},
   issn={0022-0396},
}
\bib{G}{article}{
   author={Gripenberg, Gustaf},
   title={Volterra integro-differential equations with accretive
   nonlinearity},
   journal={J. Differential Equations},
   volume={60},
   date={1985},
   number={1},
   pages={57--79},
   issn={0022-0396},
}
\bib{GLS}{article}{
   author={Gripenberg, G},
   author={Londen, S.-O},
   author={Staffans, O},
   title={Volterra integral and functional equations},
   journal={Encyclopedia of Mathematics and its Applications},
   volume={34},
    publisher={Cambridge University Press, Cambridge},
   date={1990}
}
\bib{HP}{article}{
   author={Herrero, Miguel A.},
   author={Pierre, Michel},
   title={The Cauchy problem for $u_t=\Delta u^m$ when $0<m<1$},
   journal={Trans. Amer. Math. Soc.},
   volume={291},
   date={1985},
   number={1},
   pages={145--158},
   issn={0002-9947},
}
\bib{KSVZ}{article}{
   author={Kemppainen, Jukka},
   author={Siljander, Juhana},
   author={Vergara, Vicente},
   author={Zacher, Rico},
   title={Decay estimates for time-fractional and other non-local in time
   subdiffusion equations in $\mathbb{R}^d$},
   journal={Math. Ann.},
   volume={366},
   date={2016},
   number={3-4},
   pages={941--979},
   issn={0025-5831},
}
\bib{KSZ}{article}{
   author={Kemppainen, Jukka},
   author={Siljander, Juhana},
   author={Zacher, Rico},
   title={Representation of solutions and large-time behavior for fully
   nonlocal diffusion equations},
   journal={J. Differential Equations},
   volume={263},
   date={2017},
   number={1},
   pages={149--201},
   issn={0022-0396},
}
 \bib{KY}{article}{
   author={Kubica, Adam},
   author={Yamamoto, Masahiro},
   title={Initial-boundary value problems for fractional diffusion equations
   with time-dependent coefficients},
   journal={Fract. Calc. Appl. Anal.},
   volume={21},
   date={2018},
   number={2},
   pages={276--311},
   issn={1311-0454},
}
\bib{LL}{article}{
   author={Li, Lei},
   author={Liu, Jian-Guo},
   title={A discretization of Caputo derivatives with application to time
   fractional SDEs and gradient flows},
   journal={SIAM J. Numer. Anal.},
   volume={57},
   date={2019},
   number={5},
   pages={2095--2120},
   issn={0036-1429},
}
\bib{L}{article}{
   author={Luchko, Yury},
   title={Some uniqueness and existence results for the
   initial-boundary-value problems for the generalized time-fractional
   diffusion equation},
   journal={Comput. Math. Appl.},
   volume={59},
   date={2010},
   number={5},
   pages={1766--1772},
   issn={0898-1221},
}
\bib{LY}{article}{
   author={Luchko, Yuri},
   author={Yamamoto, Masahiro},
   title={Maximum principle for the time-fractional PDEs},
   conference={
      title={Handbook of fractional calculus with applications. Vol. 2},
   },
   book={
      publisher={De Gruyter, Berlin},
   },
   date={2019},
   pages={299--325},
}
\bib{P}{book}{
   author={Pr\"uss, Jan},
   title={Evolutionary integral equations and applications},
   series={Monographs in Mathematics},
   volume={87},
   publisher={Birkh\"auser Verlag, Basel},
   date={1993},
   pages={xxvi+366},
}
\bib{PL}{article}{
   author={P\l ociniczak, \L ukasz},
   title={Analytical studies of a time-fractional porous medium equation.
   Derivation, approximation and applications},
   journal={Commun. Nonlinear Sci. Numer. Simul.},
   volume={24},
   date={2015},
   number={1-3},
   pages={169--183},
   issn={1007-5704},
}
\bib{SW}{article}{
   author={Schmitz, Kerstin},
   author={Wittbold, Petra},
   title={Entropy solutions for time-fractional porous medium type
   equations},
   journal={Differential Integral Equations},
   volume={37},
   date={2024},
   number={5-6},
   pages={309--322},
   issn={0893-4983},
}
\bib{S}{book}{
   author={Showalter, R. E.},
   title={Monotone operators in Banach space and nonlinear partial
   differential equations},
   series={Mathematical Surveys and Monographs},
   volume={49},
   publisher={American Mathematical Society, Providence, RI},
   date={1997},
   pages={xiv+278},
   isbn={0-8218-0500-2},
}
\bib{VWZ}{article}{
      title={Duality estimates for subdiffusion problems including time-fractional porous medium type equations}, 
      author={Viana, Arl\'{u}cio},
      author={Wolejko, Patryk},
      author={Zacher, Rico},
      year={2025},
      eprint={2509.07862},
}
\bib{V}{article}{
      title={A survey on mass conservation and related topics in nonlinear diffusion}, 
      author={V\'azquez, Juan Luis},
      year={2023},
      eprint={2311.18357},
}
\bib{V1}{book}{
   author={V\'azquez, Juan Luis},
   title={Smoothing and decay estimates for nonlinear diffusion equations},
   series={Oxford Lecture Series in Mathematics and its Applications},
   volume={33},
   publisher={Oxford University Press, Oxford},
   date={2006},
   pages={xiv+234},
   isbn={978-0-19-920297-3},
   isbn={0-19-920297-4},
}
\bib{V2}{book}{
   author={V\'azquez, Juan Luis},
   title={The porous medium equation},
   series={Oxford Mathematical Monographs},
   publisher={The Clarendon Press, Oxford University Press, Oxford},
   date={2007},
   pages={xxii+624},
   isbn={978-0-19-856903-9},
   isbn={0-19-856903-3},
}
\bib{VZ1}{article}{
   author={Vergara, Vicente},
   author={Zacher, Rico},
   title={Lyapunov functions and convergence to steady state for
   differential equations of fractional order},
   journal={Math. Z.},
   volume={259},
   date={2008},
   number={2},
   pages={287--309},
   issn={0025-5874},
}
\bib{VZ2}{article}{
   author={Vergara, Vicente},
   author={Zacher, Rico},
   title={Optimal decay estimates for time-fractional and other nonlocal subdiffusion equations via energy methods},
   journal={SIAM J. Math. Anal.},
   volume={47},
   date={2015},
   number={1},
   pages={210--239},
   issn={0036-1410},
}
\bib{VZ3}{article}{
   author={Vergara, Vicente},
   author={Zacher, Rico},
   title={Stability, instability, and blowup for time fractional and other
   nonlocal in time semilinear subdiffusion equations},
   journal={J. Evol. Equ.},
   volume={17},
   date={2017},
   number={1},
   pages={599--626},
   issn={1424-3199},
}
\bib{WWZ}{article}{
   author={Wittbold, Petra},
   author={Wolejko, Patryk},
   author={Zacher, Rico},
   title={Bounded weak solutions of time-fractional porous medium type and more general nonlinear and degenerate evolutionary integro-differential equations},
   journal={J. Math. Anal. Appl.},
   volume={499},
   date={2021},
   number={1},
   pages={Paper No. 125007, 20},
   issn={0022-247X},
}
\bib{Z1}{article}{
   author={Zacher, Rico},
   title={A De Giorgi--Nash type theorem for time fractional diffusion
   equations},
   journal={Math. Ann.},
   volume={356},
   date={2013},
   number={1},
   pages={99--146},
   issn={0025-5831},
}
\bib{Z2}{article}{
   author={Zacher, Rico},
   title={Time fractional diffusion equations: solution concepts,
   regularity, and long-time behavior},
   conference={
      title={Handbook of fractional calculus with applications. Vol. 2},
   },
   book={
      publisher={De Gruyter, Berlin},
   },
   isbn={978-3-11-057082-3},
   isbn={978-3-11-057166-0},
   isbn={978-3-11-057105-9},
   date={2019},
   pages={159--179},
}
\bib{Z3}{article}{
   author={Zacher, Rico},
   title={Boundedness of weak solutions to evolutionary partial
   integro-differential equations with discontinuous coefficients},
   journal={J. Math. Anal. Appl.},
   volume={348},
   date={2008},
   number={1},
   pages={137--149},
   issn={0022-247X},
}
\bib{Z4}{article}{
   author={Zacher, Rico},
   title={Weak solutions of abstract evolutionary integro-differential
   equations in Hilbert spaces},
   journal={Funkcial. Ekvac.},
   volume={52},
   date={2009},
   number={1},
   pages={1--18},
   issn={0532-8721},
}
\end{biblist}
\end{bibdiv}
\end{document}